\documentclass[12pt]{article}         \usepackage{amsfonts} \usepackage{latexsym} \usepackage{amssymb} \usepackage{graphicx} \usepackage{amsbsy}

\newtheorem{thm}{Theorem}

\newcounter{claimcount}[thm]  
\newtheorem{prop}[thm]{Proposition}

\newtheorem{cor}[thm]{Corollary}

\newcommand{\bprf}[1][Proof:]{\begin{list}{} 			{\setlength{\leftmargin}{1em} 			\setlength{\rightmargin}{0em}}                         \item {\bf \hspace{-1em}  #1 \ \ }} 
\newfont{\cursive}{cmfi10 scaled \magstep0} 

\begin{document}
\title{Conservation Laws in Cellular Automata}
\author{Marcus Pivato\\ {\em Department of Mathematics, Trent University} \\  {\em Peterborough, Ontario,  K9L 1Z6, Canada} \\
 email: {\tt mpivato@trentu.ca} or {\tt pivato@math.uh.edu}\\
 Phone: (705) 748 1011 x1293 \ \ \ \   Fax: (705)  748 1630}

\maketitle

\begin{abstract} If ${\mathbb{X}}$ is a discrete abelian group and ${\mathcal{ A}}$ a
finite set, then a cellular automaton (CA) is a continuous map
${\mathfrak{ F}}:{\mathcal{ A}}^{\mathbb{X}}{{\longrightarrow}}{\mathcal{ A}}^{\mathbb{X}}$ that commutes with all
${\mathbb{X}}$-shifts.  If $\phi:{\mathcal{ A}}{{\longrightarrow}}{\mathbb{R}}$, then, for any
${\mathbf{ a}}\in{\mathcal{ A}}^{\mathbb{X}}$, we define
$\Sigma\phi({\mathbf{ a}})=\sum_{{\mathsf{ x}}\in{\mathbb{X}}}\phi(a_{\mathsf{ x}})$ (if finite); \ $\phi$
is {\em conserved} by ${\mathfrak{ F}}$ if $\Sigma\phi$ is constant under the
action of ${\mathfrak{ F}}$.

  We characterize such {\em conservation laws} in several ways,
deriving both theoretical consequences and practical tests, and
provide a method for constructing all one-dimensional CA exhibiting a
given conservation law.
\end{abstract}
  If ${\mathcal{ A}}$ is a finite set (with discrete topology), and ${\mathbb{X}}$ an
arbitrary indexing set, then ${\mathcal{ A}}^{\mathbb{X}}$ (the space of all functions
${\mathbb{X}}\mapsto{\mathcal{ A}}$) is compact and totally disconnected in the
Tychonoff topology.  If $({\mathbb{X}},+)$ is a discrete abelian
group\footnote{We assume ${\mathbb{X}}$ is abelian only for expositional simplicity;
all results extend easily to nonabelian ${\mathbb{X}}$.} (eg. ${\mathbb{X}}={\mathbb{Z}}^D$)
 with identity ${\mathsf{ O}}$,
then ${\mathbb{X}}$ acts on itself by {\bf translation}; \ this induces a
{\bf shift} action of ${\mathbb{X}}$ on ${\mathcal{ A}}^{\mathbb{X}}$: if
${\mathbf{ a}}={\left[a_{\mathsf{ x}}  |_{{\mathsf{ x}}\in{\mathbb{X}}}^{} \right]}\in{\mathcal{ A}}^{\mathbb{X}}$, and ${\mathsf{ u}}\in{\mathbb{X}}$,
then $ {{{\boldsymbol{\sigma}}}^{{\mathsf{ u}}}} ({\mathbf{ a}}) =
{\left[b_{\mathsf{ x}}  |_{{\mathsf{ x}}\in{\mathbb{X}}}^{} \right]}$, where $b_{\mathsf{ x}} = a_{({\mathsf{ x}}+{\mathsf{ u}})}$. 

  A {\bf cellular automaton} (CA) is a continuous map
 ${\mathfrak{ F}}:{\mathcal{ A}}^{\mathbb{X}}\!\longrightarrow{\mathcal{ A}}^{\mathbb{X}}\!$ which commutes with all shifts.  The
 Curtis-Hedlund-Lyndon Theorem
\cite{HedlundCA} says that ${\mathfrak{ F}}$ is a CA if and only if there is
some finite ${\mathbb{B}}\subset{\mathbb{X}}$ (a ``neighbourhood of the
identity'') and a {\bf local map} ${\mathfrak{ f}}:{\mathcal{ A}}^{\mathbb{B}}{{\longrightarrow}}{\mathcal{ A}}$ so that,
for all ${\mathbf{ a}}\in{\mathcal{ A}}^{\mathbb{X}}$ and ${\mathsf{ x}}\in{\mathbb{X}}$, \ ${\mathfrak{ F}}({\mathbf{ a}})_{\mathsf{ x}} =
{\mathfrak{ f}}\left({\mathbf{ a}}\raisebox{-0.3em}{$\left|_{{\mathbb{B}}+{\mathsf{ x}}}\right.$}\right)$.  Here, for any ${\mathbb{W}}\subset{\mathbb{X}}$,
we define ${\mathbf{ a}}\raisebox{-0.3em}{$\left|_{{\mathbb{W}}}\right.$} =
{\left[a_{\mathsf{ w}}  |_{{\mathsf{ w}}\in{\mathbb{W}}}^{} \right]}\in{\mathcal{ A}}^{\mathbb{W}}$.  For example, if
${\mathbb{X}}={\mathbb{Z}}$ and ${\mathbb{B}}={\left[ -B...B \right]}$, then for any ${\mathsf{ z}}\in{\mathbb{Z}}$,
\ \ ${\mathbf{ a}}\raisebox{-0.3em}{$\left|_{{\mathbb{B}}+{\mathsf{ z}}}\right.$} = \left[a_{{\mathsf{ z}}-B},\ldots,a_{{\mathsf{ z}}+B}\right]$.
Without loss of generality, we assume ${\mathbb{B}}$ is symmetric, in the
sense that $\left({\mathsf{ b}}\in{\mathbb{B}}\right) \iff
\left(-{\mathsf{ b}}\in{\mathbb{B}}\right)$.

  Let $({\mathcal{ Z}},+)$ be an abelian group (usually ${\mathbb{Z}}$), and
let $\phi:{\mathcal{ A}}{{\longrightarrow}}{\mathcal{ Z}}$.  Heuristically speaking, $\phi(a)$ measures
the ``content'' of a cell in state $a\in{\mathcal{ A}}$.  If $0$ is the
identity of ${\mathcal{ Z}}$,  we refer to
$\mathbf{0}=\phi^{-1}\{0\}$ as the set of {\bf vacuum} states, 
and $\mathbf{0}^{\mathbb{X}}={\left\{ {\mathbf{ a}}\in{\mathcal{ A}}^{\mathbb{X}} \; ; \; a_{\mathsf{ x}}\in\mathbf{0}, \ \forall{\mathsf{ x}}\in{\mathbb{X}} \right\} }$
as the set of {\bf global vacuum} configurations.
If ${\mathbf{ a}}\in{\mathcal{ A}}^{\mathbb{X}}$, then the
function $\phi({\mathbf{ a}}):{\mathbb{X}}{{\longrightarrow}}{\mathcal{ Z}}$ is defined:
$\phi_{\mathsf{ x}}({\mathbf{ a}})=\phi(a_{\mathsf{ x}})$.  The {\bf support} of ${\mathbf{ a}}$ is the set
${\sf supp}\left[{\mathbf{ a}}\right]={\left\{ {\mathsf{ x}}\in{\mathbb{X}} \; ; \; a_{\mathsf{ x}}\not\in\mathbf{0} \right\} }$; \ let
${\mathcal{ A}}^{<{\mathbb{X}}}$ be the set of elements of ${\mathcal{ A}}^{\mathbb{X}}$ with finite
support.  A CA ${\mathfrak{ F}}:{\mathcal{ A}}^{\mathbb{X}}\!\longrightarrow{\mathcal{ A}}^{\mathbb{X}}\!$ is {\bf vacuum-preserving} if
${\mathfrak{ F}}(\mathbf{0}^{\mathbb{X}})\subset\mathbf{0}^{\mathbb{X}}$, or, equivalently, if
${\mathfrak{ F}}\left({\mathcal{ A}}^{<{\mathbb{X}}}\right) \subset {\mathcal{ A}}^{<{\mathbb{X}}}$.

Define $\Sigma\phi:{\mathcal{ A}}^{<{\mathbb{X}}}{{\longrightarrow}}{\mathcal{ Z}}$ by:
$\displaystyle \ \ \Sigma\phi({\mathbf{ a}}) \ = \
\sum_{{\mathsf{ x}}\in{\mathbb{X}}} \phi_{\mathsf{ x}}({\mathbf{ a}})$.
If ${\mathfrak{ F}}:{\mathcal{ A}}^{\mathbb{X}}\!\longrightarrow{\mathcal{ A}}^{\mathbb{X}}\!$ \ is a CA, then
we say $\phi$ is {\bf conserved} by ${\mathfrak{ F}}$ if, for any
 ${\mathbf{ a}}\in{\mathcal{ A}}^{<{\mathbb{X}}}$, \ $\displaystyle \Sigma\phi{\mathfrak{ F}}({\mathbf{ a}}) \ = \ \Sigma\phi({\mathbf{ a}})$; \ 
we then write: $\phi\in{\mathcal{ C}}({\mathfrak{ F}};{\mathcal{ Z}})$.
Note that ${\mathfrak{ F}}$ must be vacuum-preserving to conserve $\phi$.

        \refstepcounter{thm}                     \begin{list}{} 			{\setlength{\leftmargin}{1em} 			\setlength{\rightmargin}{0em}}                     \item {\bf Example \thethm:} Let ${\mathcal{ A}}={\mathcal{ Z}}={{\mathbb{Z}}_{/2}}$, and let $\phi:{\mathcal{ A}}{{\longrightarrow}}{\mathcal{ Z}}$ be the identity.
If ${\mathbb{X}}={\mathbb{Z}}$ and ${\mathbf{ a}}\in{\mathcal{ A}}^{<{\mathbb{Z}}}$,
then $\displaystyle \Sigma\phi({\mathbf{ a}}) = \sum_{{\mathsf{ n}}=-{\infty}}^{{\infty}} a_{\mathsf{ n}}$ measures the {\em parity}
of ${\mathbf{ a}}$ as a sequence of binary digits.
If $\mathbb{B}=\{-1,0,1\}$, \ and $\mathfrak{ f}(a_{-1},a_0,a_1) = a_{-1}+a_0+a_1$, \
then $\mathfrak{F}$ is  parity-preserving: \ 
$\Sigma\phi\circ\mathfrak{F} = \Sigma\phi$.
  	\hrulefill\end{list}   

  Examples of ${\mathbb{Z}}$-valued conservation laws for simple CA on
$\left({{\mathbb{Z}}_{/2}}\right)^{\mathbb{Z}}$ are described in \cite{KotzeSteeb}.
Necessary and sufficient conditions for conservation laws on
one-dimensional CA are given in \cite{HattoriTakesue}, and used to
completely enumerate the conservation laws for the 256 ``elementary''
(ie. nearest-neighbour) CA on $\left({{\mathbb{Z}}_{/2}}\right)^{\mathbb{Z}}$, and the 256
``elementary reversible'' CA \cite{TakesueERCA}.

  Conservation laws arise most frequently in the context of {\bf
particle-preserving} cellular automata (PPCA).  If ${\mathcal{ Z}}={\mathbb{Z}}$ and
$\phi:{\mathcal{ A}}{{\longrightarrow}}{\mathbb{N}}$, then we interpret $\phi(a_{\mathsf{ x}})$ as 
the number of ``particles'' at site ${\mathsf{ x}}\in{\mathbb{X}}$.  Thus, $\Sigma\phi$
tallies the total number of particles in space; \ $\phi$ is conserved if
particles are neither created nor
destroyed.  By extension, if ${\mathcal{ Z}}={\mathbb{Z}}^K$ and
$\phi:{\mathcal{ A}}{{\longrightarrow}}{\mathbb{N}}^K$, then $\phi$ simultaneously tallies
$K$ distinct species of indestructible particles.  In
the simplest PPCA, $\phi:{\mathcal{ A}}{{\longrightarrow}}\{0,1\}$ (e.g.
${\mathcal{ A}}=\{0,1\}$, and $\phi$ is identity map); \ thus, at most
one particle can occupy any site.  
PPCA on $\{0,1\}^{\mathbb{Z}}$ appear as models of
traffic flow
\cite{EsserSchreckenberg,NagelSchrekenberg,FukuiIshibashi,BoccaraFuks1,BoccaraFuks2,SimonNagel}, and eutectic alloys \cite{Kohyama1,Kohyama2}.

        \refstepcounter{thm}                     \begin{list}{} 			{\setlength{\leftmargin}{1em} 			\setlength{\rightmargin}{0em}}                     \item {\bf Example \thethm:} \label{EXPPCA} If ${\mathbb{X}}={\mathbb{Z}}$ and ${\mathbb{B}}=\{-1,0,1\}$, then there
  are exactly five PPCA with local maps
  ${\mathfrak{ f}}:\{0,1\}^{\mathbb{B}}{{\longrightarrow}}\{0,1\}$.  These are the identity map, the
  left- and right-shifts, and CA numbers 184 and 226 in the
  Wolfram nomenclature
\cite{WolframBook}.  In CA \#184, each ``1'' particle moves to the
right whenever there is a ``0'' to its right, and remains
stationary if there is a ``1'' to the right.  CA \#26 is the
mirror image, with movement to the left. \cite{BoccaraFuks1}.
  	\hrulefill\end{list}

  In \S\ref{S:finitary}, we characterize ${\mathcal{ C}}({\mathfrak{ F}};{\mathcal{ Z}})$ in terms
amenable to computational testing on a finite spatial domain.  In
\S\ref{S:recode}, we use this to show how any $\phi:{\mathcal{ A}}{{\longrightarrow}}{\mathbb{R}}$ can
be ``recoded'' by a ${\mathbb{R}}^+$-valued or ${\mathbb{N}}^K$-valued function
having equivalent conservation properties.  In \S\ref{S:infinitary},
we characterize ${\mathcal{ C}}({\mathfrak{ F}};{\mathbb{R}}^+)$ in terms of configurations with
infinite support.  In \S\ref{S:Cesaro} we characterize
${\mathcal{ C}}({\mathfrak{ F}};{\mathbb{R}}^+)$ in terms of spatial ergodic averages, assuming
${\mathbb{X}}$ is an amenable group, while in \S\ref{S:measure}, we
characterize ${\mathcal{ C}}({\mathfrak{ F}})$ in terms of stationary measures on
${\mathcal{ A}}^{\mathbb{X}}$, even when ${\mathbb{X}}$ is not amenable.  In
\S\ref{S:construct}, we consider the construction of CA with a
particular conservation law.

\section{Finitary Characterizations
\label{S:finitary}}

  Let $\mbox{\cursive O}\in\mathbf{0}$ be some fixed vacuum state.  If ${\mathbb{V}}\subset{\mathbb{X}}$
 is finite and ${\mathbf{ a}}\in{\mathcal{ A}}^{\mathbb{V}}$, then let ${\left\langle {\mathbf{ a}} \right\rangle }$ denote the
 configuration ${\mathbf{ b}}\in{\mathcal{ A}}^{<{\mathbb{X}}}$ defined by: ${\mathbf{ b}}\raisebox{-0.3em}{$\left|_{{\mathbb{V}}}\right.$}={\mathbf{ a}}$
 and $b_{\mathsf{ x}}=\mbox{\cursive O}$ for all ${\mathsf{ x}}\not\in{\mathbb{V}}$.  Then define
 $\displaystyle\Sigma\phi({\mathbf{ a}}) = \Sigma\phi{\left\langle {\mathbf{ a}} \right\rangle }$.  If ${\mathbf{ a}}\in{\mathcal{ A}}^{\mathbb{X}}$,
 then define $\displaystyle \Sigma\phi({\mathbf{ a}})\raisebox{-0.3em}{$\left|_{{\mathbb{V}}}\right.$} = \sum_{{\mathsf{ v}}\in{\mathbb{V}}}
\phi\left(a_{\mathsf{ v}}\right)$.  If ${\mathbb{U}},{\mathbb{V}}\subset{\mathbb{X}}$, define
 ${\mathbb{V}}+{\mathbb{U}}={\left\{ {\mathsf{ v}}+{\mathsf{ u}} \; ; \; {\mathsf{ v}} \in{\mathbb{V}}, \ {\mathsf{ u}}\in{\mathbb{U}} \right\} }$.  If
 ${\mathbf{ a}}\in{\mathcal{ A}}^{{\mathbb{V}}+{\mathbb{B}}}$, then we define ${\mathfrak{ F}}({\mathbf{ a}})\in{\mathcal{ A}}^{\mathbb{V}}$ by:
 ${\mathfrak{ F}}({\mathbf{ a}})_{\mathsf{ v}} = {\mathfrak{ f}}\left({\mathbf{ a}}\raisebox{-0.3em}{$\left|_{{\mathsf{ v}}+{\mathbb{B}}}\right.$}\right)$ for all
 ${\mathsf{ v}}\in{\mathbb{V}}$. 
Let ${{\mathbb{B}}^{(2)}} = {\mathbb{B}} + {{\mathbb{B}}}$; \  for example, if
${\mathbb{B}}={\left[ -B_0...B_1 \right]}^D\subset{\mathbb{Z}}^D$, then ${{\mathbb{B}}^{(2)}}={\left[ -2B_0...2B_1 \right]}^D$.

\begin{prop}\label{finitary.thm}  
 $\phi\in{\mathcal{ C}}({\mathfrak{ F}})$ if and only if, for all ${\mathbf{ a}},{\mathbf{ c}} \in {\mathcal{ A}}^{{\mathbb{B}}^{(2)}}$,
identical everywhere except that  $a_{\mathsf{ o}} = a\neq c=c_{\mathsf{ o}}$,
we have:
  \begin{equation} \label{foobar}
\Sigma\phi \ {\mathfrak{ F}}({\mathbf{ c}}) - \Sigma\phi\ {\mathfrak{ F}}({\mathbf{ a}})  
 \ \ = \ \ 
\phi(c) - \phi(a).
 \end{equation}
 \end{prop}
\bprf
Let $\delta = \phi(c) - \phi(a)$.

{\bf \hspace{-1em}  Proof of ``$\mbox{$\Longrightarrow$}$'': \ \ }  Clearly, $\Sigma\phi\left({\left\langle {\mathbf{ c}} \right\rangle }\right) =
 \delta+ \Sigma\phi\left({\left\langle {\mathbf{ a}} \right\rangle }\right)$, and thus,
 $\Sigma\phi{\mathfrak{ F}}\left({\left\langle {\mathbf{ c}} \right\rangle }\right)\ =
 \delta+\Sigma\phi\,{\mathfrak{ F}}\left({\left\langle {\mathbf{ a}} \right\rangle }\right)$.  Now, ${\left\langle {\mathbf{ a}} \right\rangle }$ and
 ${\left\langle {\mathbf{ c}} \right\rangle }$ only differ at ${\mathsf{ O}}$, so 
${\mathfrak{ F}}\left({\left\langle {\mathbf{ a}} \right\rangle }\right)\raisebox{-0.3em}{$\left|_{{\mathbb{X}}\setminus{\mathbb{B}}}\right.$}={\mathfrak{ F}}\left({\left\langle {\mathbf{ c}} \right\rangle }\right)\raisebox{-0.3em}{$\left|_{{\mathbb{X}}\setminus{\mathbb{B}}}\right.$}$, while  
${\mathfrak{ F}}\left({\left\langle {\mathbf{ a}} \right\rangle }\right)\raisebox{-0.3em}{$\left|_{{\mathbb{B}}}\right.$}={\mathfrak{ F}}\left({\mathbf{ a}}\right)$ \ and \ 
${\mathfrak{ F}}\left({\left\langle {\mathbf{ c}} \right\rangle }\right)\raisebox{-0.3em}{$\left|_{{\mathbb{B}}}\right.$}={\mathfrak{ F}}\left({\mathbf{ c}}\right)$.  Thus, \ 
\[ 
\Sigma \phi \ {\mathfrak{ F}}\left({\mathbf{ c}}\right)  -  \Sigma \phi\ {\mathfrak{ F}}\left({\mathbf{ a}}\right) 
\quad = \quad
\Sigma \phi \ {\mathfrak{ F}}\left({\left\langle {\mathbf{ c}} \right\rangle }\right)\raisebox{-0.3em}{$\left|_{{\mathbb{B}}}\right.$} 
   -\Sigma \phi\ {\mathfrak{ F}}\left({\left\langle {\mathbf{ a}} \right\rangle }\right)\raisebox{-0.3em}{$\left|_{{\mathbb{B}}}\right.$}
 \quad = \quad 
\Sigma \phi\ {\mathfrak{ F}}\left({\left\langle {\mathbf{ c}} \right\rangle }\right)  -  \Sigma \phi\ {\mathfrak{ F}}\left({\left\langle {\mathbf{ a}} \right\rangle }\right)
 \ = \quad \delta,
\]
  which yields equation (\ref{foobar}). 

{\bf \hspace{-1em}  Proof of ``${\Longleftarrow}$'': \ \ }  Suppose ${\mathbf{ a}}\in{\mathcal{ A}}^{<{\mathbb{X}}}$,
with $\displaystyle{\sf supp}\left[{\mathbf{ a}}\right]=\{{\mathsf{ x}}_1,{\mathsf{ x}}_2,\ldots,{\mathsf{ x}}_N\}$.  For all $n\in{\left[ 1..N \right]}$,
suppose $a_{{\mathsf{ x}}_n} = b_n$ and let $\delta_n = \phi(b_n)$.

Consider the vacuum ${\mathbf{ a}}^{0} \in\mathbf{0}^{\mathbb{X}}$ defined:
\[
  a^{0}_{\mathsf{ x}} \ =
{\left\{ \begin{array}{rcl}                                  \mbox{\cursive O}  &&  \ \ \mbox{if} \ {\mathsf{ x}}\in{\sf supp}\left[{\mathbf{ a}}\right]\\
	 a_{\mathsf{ x}} && \mbox{otherwise.}\\                                \end{array}  \right.  }
\]
 We build
${\mathbf{ a}}$ from ${\mathbf{ a}}^{0}$ one nonvacuum site at a time.  For 
$n\in{\left[ 1..N \right]}$, define \ ${\mathbf{ a}}^{n}\in{\mathcal{ A}}^{<{\mathbb{X}}}$ by:
\[
  a^{n}_{\mathsf{ x}} \ =
{\left\{ \begin{array}{rcl}                                   b_k &&  \ \ \mbox{if} \ {\mathsf{ x}}={\mathsf{ x}}_k \ \mbox{for some} \ k\leq n,\\
	 a^{0}_{\mathsf{ x}} && \mbox{otherwise.}\\                                \end{array}  \right.  }
\]
 Thus, ${\mathbf{ a}}^{N} = {\mathbf{ a}}$.
 For any $n\in{\left[ 0..N \right)}$, \ 
${\mathbf{ a}}^{n}$ and ${\mathbf{ a}}^{n+1}$ differ only at ${\mathsf{ x}}_n$, so
${\mathfrak{ F}}\left({\mathbf{ a}}^{n}\right)$ and ${\mathfrak{ F}}\left({\mathbf{ a}}^{n+1}\right)$ differ only 
in ${\mathsf{ x}}_n+{{\mathbb{B}}}$; \ hence
\begin{eqnarray*}
\Sigma\phi\ {\mathfrak{ F}}\left({\mathbf{ a}}^{n}\right)
\ - 
\Sigma\phi\ {\mathfrak{ F}}\left({\mathbf{ a}}^{n-1}\right)
& = &
\Sigma\phi\ {\mathfrak{ F}}\left({\mathbf{ a}}^{n}\right)\raisebox{-0.3em}{$\left|_{{\mathsf{ x}}_n+{{\mathbb{B}}}}\right.$}
\ - 
\Sigma\phi\ {\mathfrak{ F}}\left({\mathbf{ a}}^{n-1}\right)\raisebox{-0.3em}{$\left|_{{\mathsf{ x}}_n+{{\mathbb{B}}}}\right.$}\\
&= & \Sigma\phi\ {\mathfrak{ F}}\left({\mathbf{ a}}^{n}\raisebox{-0.3em}{$\left|_{{\mathsf{ x}}_n+{{\mathbb{B}}^{(2)}}}\right.$}\right)
\ - 
\Sigma\phi\ {\mathfrak{ F}}\left({\mathbf{ a}}^{n-1}\raisebox{-0.3em}{$\left|_{{\mathsf{ x}}_n+{{\mathbb{B}}^{(2)}}}\right.$}\right)\\
&=_{[1]} &
 \phi\ \left(a^{n}_{{\mathsf{ x}}_n}\right)
\ - \
\phi\ \left(a^{n-1}_{{\mathsf{ x}}_n}\right)
\ \   =  \ \  \delta_n,
\end{eqnarray*}
where [1] follows from applying equation (\ref{foobar}) at ${\mathsf{ x}}_n$.
Inductively,  $\displaystyle \Sigma\phi\ {\mathfrak{ F}}\left({\mathbf{ a}}^{N}\right)  \ = \
 \sum_{n=1}^N \delta_n  \ + \Sigma\phi\ {\mathfrak{ F}}\left({\mathbf{ a}}^{0}\right) \ =  \
 \sum_{n=1}^N \delta_n \ = \  \Sigma\phi\left[{\mathbf{ a}}\right]$.
 {\tt \hrulefill $\Box$ } \end{list}  \medskip  

Proposition \ref{finitary.thm} generalizes Proposition 2.3 of
\cite{HattoriTakesue} (which is the case  ${\mathbb{X}}={\mathbb{Z}}$).  
There is also a characterization of ${\mathbb{N}}$-valued conservation laws
in terms of periodic configurations (see Theorem 2.1 of
\cite{BoccaraFuks2} for case ${\mathbb{X}}={\mathbb{Z}}$ or Proposition 1 of
\cite{DurandFormentiRoka} for case ${\mathbb{X}}={\mathbb{Z}}^D$); \ we generalize
this to the following characterization of arbitrary conservation laws
for any group ${\mathbb{X}}$.

   Let ${\mathfrak{ Q}}$ be the set of all finite quotient groups ${\widetilde{\mathbb{X}}}$ of
${\mathbb{X}}$ such that ${{\mathbb{B}}^{(2)}}$ maps bijectively onto its image
${{\widetilde{\mathbb{B}}}^{(2)}}\subset{\widetilde{\mathbb{X}}}$ under the quotient map.  The local map
${\mathfrak{ f}}:{\mathcal{ A}}^{\mathbb{B}}{{\longrightarrow}}{\mathcal{ A}}$ induces a cellular automaton
${\widetilde{\mathfrak{ F}}}:{\mathcal{ A}}^{\widetilde{\mathbb{X}}}\!\longrightarrow{\mathcal{ A}}^{\widetilde{\mathbb{X}}}\!$ for any ${\widetilde{\mathbb{X}}}\in{\mathfrak{ Q}}$.  Let
${\mathcal{ C}}\left({\widetilde{\mathfrak{ F}}}\right)\ =\ {\left\{ \phi:{\mathcal{ A}}{{\longrightarrow}}{\mathbb{R}} \; ; \; \phi\ \mbox{is
conserved by}\ {\widetilde{\mathfrak{ F}}} \right\} }.$

\begin{cor}\label{fuks.boccara}  
  $\phi\in{\mathcal{ C}}({\mathfrak{ F}})$ iff $\phi\in{\mathcal{ C}}\left({\widetilde{\mathfrak{ F}}}\right)$.
 \end{cor}
\bprf 
By Proposition \ref{finitary.thm}, $\phi\in{\mathcal{ C}}({\mathfrak{ F}})$ iff equation
(\ref{foobar}) is true, while and $\phi\in{\mathcal{ C}}\left({\widetilde{\mathfrak{ F}}}\right)$ iff
$\widetilde{(\ref{foobar})}$ is true, where
$\widetilde{(\ref{foobar})}$ is (\ref{foobar}) with ``${{\widetilde{\mathbb{B}}}^{(2)}}$''
replacing ``${{\mathbb{B}}^{(2)}}$''.  If we identify ${\mathcal{ A}}^{{{\widetilde{\mathbb{B}}}^{(2)}}}$ and
${\mathcal{ A}}^{{\mathbb{B}}^{(2)}}$ in the obvious way, then clearly,
$\widetilde{(\ref{foobar})}$ and (\ref{foobar}) are equivalent.
 {\tt \hrulefill $\Box$ } \end{list}  \medskip  

   For example, if ${\mathbb{X}}={\mathbb{Z}}^D$, and ${\mathbb{B}}={\left[ -B...B \right]}^D$, then
${\mathfrak{ Q}}$ includes the quotient group ${\widetilde{\mathbb{X}}}=({{\mathbb{Z}}_{/M}})^D$ for any
$M > 4\cdot B$.  Elements of ${\mathcal{ A}}^{\widetilde{\mathbb{X}}}$ correspond to $M$-periodic
configurations in ${\mathcal{ A}}^{{\mathbb{Z}}^D}$ (where ${\mathbf{ a}}\in{\mathcal{ A}}^{{\mathbb{Z}}^D}$ is {\bf
$M$-periodic} if $ {{{\boldsymbol{\sigma}}}^{M\cdot{\mathsf{ z}}}} ({\mathbf{ a}})={\mathbf{ a}}$ for any
${\mathsf{ z}}\in{\mathbb{Z}}^D$).  The action of ${\widetilde{\mathfrak{ F}}}$ on ${\mathcal{ A}}^{\widetilde{\mathbb{X}}}$ corresponds
to the action of ${\mathfrak{ F}}$ on $M$-periodic configurations in
${\mathcal{ A}}^{{\mathbb{Z}}^D}$.  Thus, ${\mathfrak{ F}}$ conserves $\phi$ iff ${\mathfrak{ F}}$ conserves
$\phi$ on $M$-periodic configurations.

\section{Recoding
\label{S:recode}}

  Proposition \ref{finitary.thm} yields a convenient 
``recoding'' of real-valued conservation laws.
Let ${\mathbb{R}}^+={\left\{ x\in{\mathbb{R}} \; ; \; x\geq0 \right\} }$, and let ${\mathcal{ C}}\left({\mathfrak{ F}};{\mathbb{R}}^+\right)$
denote ${\mathbb{R}}^+$-valued elements of ${\mathcal{ C}}\left({\mathfrak{ F}};{\mathbb{R}}\right)$.

\begin{prop}
\label{recode.thm}  
 Let $\phi:{\mathcal{ A}}{{\longrightarrow}}{\mathbb{R}}$.
\begin{enumerate}
\item  There is a function $\widetilde\phi:{\mathcal{ A}}{{\longrightarrow}}{\mathbb{R}}^+$ so 
that, for any cellular automaton ${\mathfrak{ F}}$,
\[\left( \ \rule[-0.5em]{0em}{1em}       \begin{minipage}{40em}       \begin{tabbing}         $\phi\in{\mathcal{ C}}({\mathfrak{ F}};{\mathbb{R}})$        \end{tabbing}      \end{minipage} \ \right) \iff
\left( \ \rule[-0.5em]{0em}{1em}       \begin{minipage}{40em}       \begin{tabbing}         $ {\widetilde{\phi }}\in{\mathcal{ C}}({\mathfrak{ F}};{\mathbb{R}}^+)$        \end{tabbing}      \end{minipage} \ \right)\]

\item  There is a function $\hat\phi:{\mathcal{ A}}{{\longrightarrow}}{\mathbb{N}}^K$ so 
that, for any cellular automaton ${\mathfrak{ F}}$,
\[
\left( \ \rule[-0.5em]{0em}{1em}       \begin{minipage}{40em}       \begin{tabbing}         $\phi\in{\mathcal{ C}}({\mathfrak{ F}};{\mathbb{R}})$        \end{tabbing}      \end{minipage} \ \right) \iff 
\left( \ \rule[-0.5em]{0em}{1em}       \begin{minipage}{40em}       \begin{tabbing}         $\hat\phi\in{\mathcal{ C}}({\mathfrak{ F}};{\mathbb{N}}^K)$        \end{tabbing}      \end{minipage} \ \right)\]
\end{enumerate}
 \end{prop}
\bprf
{\bf Part 1}: Let $\displaystyle -M=\min_{a\in{\mathcal{ A}}} \phi(a)$, and define ${\widetilde{\phi }}(a)
= \phi(a)+M$ for all $a\in{\mathcal{ A}}$.  Then clearly, ${\widetilde{\phi }}$ satisfies
the condition of Proposition \ref{finitary.thm} if and only if
$\phi$ does.

{\bf Part 2}:  ${\mathcal{ A}}$ is finite, so 
$\phi({\mathcal{ A}})\subset{\mathbb{R}}$ is finite, so  the subgroup ${\mathbb{A}}
\subset{\mathbb{R}}$ generated by $\phi({\mathcal{ A}})$ is a finitely generated,
torsion-free abelian group, therefore isomorphic to ${\mathbb{Z}}^K$ for some
$K$.  If $\zeta:{\mathbb{A}}{{\longrightarrow}}{\mathbb{Z}}^K$ is this
isomorphism, and $\hat\phi=\zeta\circ\phi$, then clearly
$\left( \ \rule[-0.5em]{0em}{1em}       \begin{minipage}{40em}       \begin{tabbing}         $\phi\in{\mathcal{ C}}({\mathfrak{ F}};{\mathbb{R}}$        \end{tabbing}      \end{minipage} \ \right)\iff
\left( \ \rule[-0.5em]{0em}{1em}       \begin{minipage}{40em}       \begin{tabbing}         $\hat\phi\in{\mathcal{ C}}({\mathfrak{ F}};{\mathbb{Z}}^K)$        \end{tabbing}      \end{minipage} \ \right)$.

  We can always  choose $\zeta$ so that $\zeta\left({\mathbb{A}}\cap{\mathbb{R}}^+\right)
\subset {\mathbb{N}}^K$, and by {\bf Part 1} we can assume $\phi$ is nonnegative,
so that $\hat\phi:{\mathcal{ A}}{{\longrightarrow}}{\mathbb{N}}^K$.
 {\tt \hrulefill $\Box$ } \end{list}  \medskip  

  Note that the vacuum states of ${\widetilde{\phi }}$ are not the same as
those of $\phi$, because  ${\widetilde{\phi }}^{-1}\{0\} = \phi^{-1}\{-M\}$.
Thus, $\phi$ and ${\widetilde{\phi }}$ determine two different notion of ``finite support'',
and ``${\mathcal{ A}}^{<{\mathbb{X}}}$'' refers to two different subsets of ${\mathcal{ A}}^{\mathbb{X}}$.

 {\bf Part 1} of Proposition \ref{recode.thm} implies that, to characterize real
conservation laws, it is sufficient to characterize nonnegative ones; \
this will be useful in \S\ref{S:infinitary} and \S\ref{S:Cesaro}.
{\bf Part 2} of Proposition \ref{recode.thm} implies that we can
interpret any real conserved quantity as tallying
$K$ species of indestructible particles.
Conversely, to {\em construct} a CA with a given real-valued
conservation law, it is sufficient to construct a $K$-species
PPCA; \ this will be useful in \S\ref{S:construct}.

\section{A Nonfinitary Characterization
\label{S:infinitary}}

   Defining conservation laws in the context of ${\mathcal{ A}}^{<{\mathbb{X}}}$
is somewhat unnatural, because ${\mathcal{ A}}^{<{\mathbb{X}}}$ is a very small subset
of ${\mathcal{ A}}^{\mathbb{X}}$.  We now characterize
conservation laws in a way which is meaningful for any ${\mathbf{ a}}\in{\mathcal{ A}}^{\mathbb{X}}$.

   For any ${\mathbb{W}}\subset{\mathbb{X}}$, define ${{\mathbf{ c}}{\mathbf{ l}}\left[ {\mathbb{W}} \right]} =
{\mathbb{B}}+{\mathbb{W}}$, and ${{\mathbf{ i}}{\mathbf{ n}}{\mathbf{ t}}\left[ {\mathbb{W}} \right]} =
{\left\{ {\mathsf{ w}}\in{\mathbb{W}} \; ; \; {\mathbb{B}}+{\mathsf{ w}}\subset{\mathbb{W}} \right\} }$.  Thus,
${\mathbb{W}}\subset{{\mathbf{ i}}{\mathbf{ n}}{\mathbf{ t}}\left[ {{\mathbf{ c}}{\mathbf{ l}}\left[ {\mathbb{W}} \right]} \right]}$.  Because ${\mathbb{B}}$ is
{symmetric}, ${{\mathbf{ i}}{\mathbf{ n}}{\mathbf{ t}}\left[ {\mathbb{W}} \right]} = {\mathbb{W}}\setminus{{\mathbf{ c}}{\mathbf{ l}}\left[ {\mathbb{W}}^\complement \right]}$, where
${\mathbb{W}}^\complement={\mathbb{X}}\setminus{\mathbb{W}}$.

\begin{thm}
\label{inf.cond}  
  Suppose $\phi:{\mathcal{ A}}{{\longrightarrow}}{\mathbb{R}}^+$.  Then $\phi\in{\mathcal{ C}}({\mathfrak{ F}})$ iff,
for any ${\mathbf{ a}}\in{\mathcal{ A}}^{\mathbb{X}}$ with ${\mathbf{ a}}'={\mathfrak{ F}}({\mathbf{ a}})$, and any finite ${\mathbb{W}}\subset{\mathbb{X}}$
\begin{equation}
\label{slow.movement}
  \Sigma \phi ({\mathbf{ a}})\raisebox{-0.3em}{$\left|_{{{\mathbf{ i}}{\mathbf{ n}}{\mathbf{ t}}\left[ {\mathbb{W}} \right]}}\right.$}
\quad  \leq \quad  
\Sigma \phi({\mathbf{ a}}')\raisebox{-0.3em}{$\left|_{{\mathbb{W}}}\right.$}
\quad  \leq \quad  
  \Sigma \phi ({\mathbf{ a}})\raisebox{-0.3em}{$\left|_{{{\mathbf{ c}}{\mathbf{ l}}\left[ {\mathbb{W}} \right]}}\right.$}
\end{equation}
 \end{thm}
\bprf
``$\mbox{$\Longrightarrow$}$'':  Let ${\mathbb{V}}={{\mathbf{ c}}{\mathbf{ l}}\left[ {\mathbb{W}} \right]}$ and
let ${\mathbf{ b}}={\left\langle {\mathbf{ a}}\raisebox{-0.3em}{$\left|_{{\mathbb{V}}}\right.$} \right\rangle }$ ---that is,
 ${\mathbf{ b}}\raisebox{-0.3em}{$\left|_{{\mathbb{V}}}\right.$} = \
{\mathbf{ a}}\raisebox{-0.3em}{$\left|_{{\mathbb{V}}}\right.$}$, and $b_{\mathsf{ x}} = \mbox{\cursive O} \in \mathbf{0}$ for
${\mathsf{ x}}\not\in{\mathbb{V}}$, so ${\mathbf{ b}}\in{\mathcal{ A}}^{<{\mathbb{X}}}$.  Thus, if ${\mathbf{ b}}'={\mathfrak{ F}}({\mathbf{ b}})$, then
${\mathbf{ b}}'\raisebox{-0.3em}{$\left|_{{\mathbb{W}}}\right.$} = {\mathbf{ a}}'\raisebox{-0.3em}{$\left|_{{\mathbb{W}}}\right.$}$, so it is sufficient to prove
(\ref{slow.movement}) for ${\mathbf{ b}}$. The right-hand inequality in
(\ref{slow.movement}) follows because:
\[
\Sigma \phi({\mathbf{ a}}')\raisebox{-0.3em}{$\left|_{{\mathbb{W}}}\right.$} \quad = \quad 
\Sigma \phi({\mathbf{ b}}')\raisebox{-0.3em}{$\left|_{{\mathbb{W}}}\right.$} \quad \leq_{_{[1]}} \quad
\Sigma \phi({\mathbf{ b}}') \quad =_{_{[2]}} \quad \Sigma \phi ({\mathbf{ b}})
\quad = \quad \Sigma \phi ({\mathbf{ a}})\raisebox{-0.3em}{$\left|_{{\mathbb{V}}}\right.$},
\]
 where $[1]$ is because $\phi$ is nonnegative, and
 $[2]$ is because $\phi\in{\mathcal{ C}}({\mathfrak{ F}})$. 

\begin{figure}
\centerline{
\includegraphics[scale=0.5]{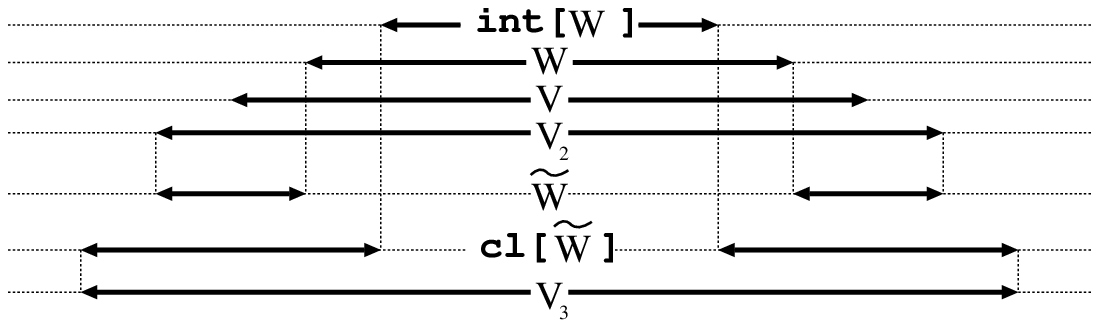}}
\caption{\label{Wfig}}
\end{figure}

 To see the left-hand inequality  in (\ref{slow.movement}),
let ${\mathbb{V}}_2={{\mathbf{ c}}{\mathbf{ l}}\left[ {\mathbb{V}} \right]}$
and ${\widetilde{\mathbb{W}}}={\mathbb{V}}_2\setminus{\mathbb{W}}$
(see Figure \ref{Wfig}). 
Thus, ${\sf supp}\left[{\mathbf{ b}}'\right]\subset{\mathbb{V}}_2$, so
\begin{equation}
\label{bb1}
\Sigma\phi ({\mathbf{ b}}') 
\quad = \quad  \Sigma\phi({\mathbf{ b}}')\raisebox{-0.3em}{$\left|_{{\mathbb{V}}_2}\right.$}
\quad = \quad  \Sigma\phi({\mathbf{ b}}')\raisebox{-0.3em}{$\left|_{{\widetilde{\mathbb{W}}}}\right.$} + \Sigma\phi({\mathbf{ b}}')\raisebox{-0.3em}{$\left|_{{\mathbb{W}}}\right.$},
\end{equation}
If ${\mathbb{V}}_3={{\mathbf{ c}}{\mathbf{ l}}\left[ {\mathbb{V}}_2 \right]}$,
then ${{\mathbf{ c}}{\mathbf{ l}}\left[ {\widetilde{\mathbb{W}}} \right]} = {\mathbb{V}}_3\setminus{{\mathbf{ i}}{\mathbf{ n}}{\mathbf{ t}}\left[ {\mathbb{W}} \right]}$, and
clearly, ${\sf supp}\left[{\mathbf{ b}}\right]\subset{\mathbb{W}}\subset{\mathbb{V}}_3$, so that
\begin{equation}
\label{bb2}
\Sigma\phi ({\mathbf{ b}}) 
\quad = \quad  \Sigma\phi({\mathbf{ b}})\raisebox{-0.3em}{$\left|_{{\mathbb{V}}_3}\right.$}
\quad = \quad  \Sigma\phi({\mathbf{ b}})\raisebox{-0.3em}{$\left|_{{{\mathbf{ c}}{\mathbf{ l}}\left[ {\widetilde{\mathbb{W}}} \right]}}\right.$} + 
\Sigma\phi({\mathbf{ b}})\raisebox{-0.3em}{$\left|_{{{\mathbf{ i}}{\mathbf{ n}}{\mathbf{ t}}\left[ {\mathbb{W}} \right]}}\right.$}.
\end{equation}
But applying the right-hand inequality  in (\ref{slow.movement})
to ${\widetilde{\mathbb{W}}}$, we have
\begin{equation}
\label{bb3}
 \Sigma\phi({\mathbf{ b}}')\raisebox{-0.3em}{$\left|_{{\widetilde{\mathbb{W}}}}\right.$} \quad \leq \quad \Sigma\phi({\mathbf{ b}})\raisebox{-0.3em}{$\left|_{{{\mathbf{ c}}{\mathbf{ l}}\left[ {\widetilde{\mathbb{W}}} \right]}}\right.$}
\end{equation}
while, by hypothesis that $\phi$ is conserved, we have
\begin{equation}
\label{bb4}
\Sigma\phi ({\mathbf{ b}})  \quad = \quad  \Sigma\phi ({\mathbf{ b}}') 
\end{equation}
 Combining (\ref{bb1}-\ref{bb4}) yields:
\begin{eqnarray*}
\Sigma\phi({\mathbf{ b}})\raisebox{-0.3em}{$\left|_{{{\mathbf{ c}}{\mathbf{ l}}\left[ {\widetilde{\mathbb{W}}} \right]}}\right.$} \ + \ 
        \Sigma\phi({\mathbf{ b}})\raisebox{-0.3em}{$\left|_{{{\mathbf{ i}}{\mathbf{ n}}{\mathbf{ t}}\left[ {\mathbb{W}} \right]}}\right.$}
&= & 
\Sigma\phi({\mathbf{ b}}')\raisebox{-0.3em}{$\left|_{{\widetilde{\mathbb{W}}}}\right.$}  \ +  \ \Sigma\phi({\mathbf{ b}}')\raisebox{-0.3em}{$\left|_{{\mathbb{W}}}\right.$}\\
& \leq &
\Sigma\phi({\mathbf{ b}})\raisebox{-0.3em}{$\left|_{{{\mathbf{ c}}{\mathbf{ l}}\left[ {\widetilde{\mathbb{W}}} \right]}}\right.$} \ +  \ \Sigma\phi({\mathbf{ b}}')\raisebox{-0.3em}{$\left|_{{\mathbb{W}}}\right.$}
\end{eqnarray*}
 from which we conclude that $ \Sigma\phi({\mathbf{ b}})\raisebox{-0.3em}{$\left|_{{{\mathbf{ i}}{\mathbf{ n}}{\mathbf{ t}}\left[ {\mathbb{W}} \right]}}\right.$}
\ \leq \
 \Sigma\phi({\mathbf{ b}}')\raisebox{-0.3em}{$\left|_{{\mathbb{W}}}\right.$}$,  as desired. 

\paragraph*{``${\Longleftarrow}$'':}
First, note that ${\mathfrak{ F}}$ must be vacuum-preserving:  If
${\mathbf{ a}}\in\mathbf{0}^{\mathbb{X}}$ and ${\mathbf{ a}}'={\mathfrak{ F}}({\mathbf{ a}})$, then
(\ref{slow.movement}) implies that, for any ${\mathsf{ x}}\in{\mathbb{X}}$,
$0\leq\phi(a'_{\mathsf{ x}}) \leq \Sigma\phi({\mathbf{ a}})\raisebox{-0.3em}{$\left|_{{{\mathbf{ c}}{\mathbf{ l}}\left[ {\mathsf{ x}} \right]}}\right.$}
= 0$, so $a'_{\mathsf{ x}}\in\mathbf{0}$.

 Next, suppose ${\mathbf{ a}}\in{\mathcal{ A}}^{<{\mathbb{X}}}$, with ${\sf supp}\left[{\mathbf{ a}}\right]={\mathbb{Y}}$.  Let
${\mathbb{W}}={{\mathbf{ c}}{\mathbf{ l}}\left[ {\mathbb{Y}} \right]}$;  \  then, since ${\mathbb{Y}}\subset{{\mathbf{ i}}{\mathbf{ n}}{\mathbf{ t}}\left[ {\mathbb{W}} \right]}$,
(\ref{slow.movement}) implies:
\[
\Sigma \phi ({\mathbf{ a}})
\quad = \quad 
  \Sigma \phi ({\mathbf{ a}})\raisebox{-0.3em}{$\left|_{{\mathbb{Y}}}\right.$}
\quad \leq \quad
  \Sigma \phi ({\mathbf{ a}})\raisebox{-0.3em}{$\left|_{{{\mathbf{ i}}{\mathbf{ n}}{\mathbf{ t}}\left[ {\mathbb{W}} \right]}}\right.$}
\quad \leq \quad
\Sigma \phi({\mathbf{ a}}')\raisebox{-0.3em}{$\left|_{{\mathbb{W}}}\right.$}
\quad \leq \quad 
  \Sigma \phi ({\mathbf{ a}})\raisebox{-0.3em}{$\left|_{{{\mathbf{ c}}{\mathbf{ l}}\left[ {\mathbb{W}} \right]}}\right.$}
\quad =  \quad 
\Sigma \phi ({\mathbf{ a}}),
\]
so that $\Sigma \phi({\mathbf{ a}}')\raisebox{-0.3em}{$\left|_{{\mathbb{W}}}\right.$} \ = \ \Sigma \phi ({\mathbf{ a}})$.
But ${\mathfrak{ F}}$ is vacuum-preserving, so ${\sf supp}\left[{\mathbf{ a}}'\right]\subset{\mathbb{W}}$; \ thus,
$\Sigma \phi({\mathbf{ a}}')\raisebox{-0.3em}{$\left|_{{\mathbb{W}}}\right.$}
\ = \ \Sigma \phi({\mathbf{ a}}')$, \ thus, \ 
$\Sigma \phi({\mathbf{ a}}') \ = \ \Sigma \phi ({\mathbf{ a}})$.  
 {\tt \hrulefill $\Box$ } \end{list}  \medskip  

\section{Conservation and Spatial Ergodic Averages
\label{S:Cesaro}}

  A {\bf F{\o}lner  sequence} \cite{Tempelman} on ${\mathbb{X}}$ is a sequence
of finite subsets ${\mathbb{I}}_n\subset{\mathbb{X}}$ so that, for any ${\mathsf{ x}}\in{\mathbb{X}}$
\[
\lim_{n{\rightarrow}{\infty}} \frac{{{\sf card}\left[{\mathbb{I}}_n\cap({\mathbb{I}}_n+{\mathsf{ x}})\right]}}{{{\sf card}\left[{\mathbb{I}}_n\right]}} = 1
\]
The group ${\mathbb{X}}$ is called {\bf amenable} if it has a F{\o}lner  sequence.
For example, ${\mathbb{Z}}^D$ is amenable, because ${\mathbb{I}}_n={\left[ 1..n \right]}^D$ 
forms a F{\o}lner  sequence.
If ${\mathbb{X}}$ is amenable and ${\mathbb{Y}}\subset{\mathbb{X}}$,
 we define the {\bf Ces\`aro  density} of 
${\mathbb{Y}}$ by
\begin{equation}
\label{Cesaro.Density}
 {{\sf density}\left[{\mathbb{Y}}\right] } \ = \ \lim_{n{\rightarrow}{\infty}}
 \frac{{{\sf card}\left[{\mathbb{Y}}\cap{\mathbb{I}}_n\right]}}{{{\sf card}\left[{\mathbb{I}}_n\right]}}.
\end{equation}
 If ${\mathbf{ a}}\in{\mathcal{ A}}^{\mathbb{X}}$ then the (spatial)
{\bf ergodic average} of $\phi$ on ${\mathbf{ a}}$ is defined:
\begin{equation}
\label{Ergodic.Average}
 \displaystyle{{\sf ErgAve}}_{{\mathbb{X}}}\  \phi({\mathbf{ a}}) \ = \ 
\lim_{n{\rightarrow}{\infty}}
 \frac{1}{I_n} \,\sum_{\ensuremath{\mathsf{ i}}\in{\mathbb{I}}_n} \phi(a_\ensuremath{\mathsf{ i}})
\hspace{3em} \mbox{(where $I_n={{\sf card}\left[{\mathbb{I}}_n\right]}$)}.
\end{equation}
If ${\mathbb{X}}={\mathbb{Z}}$ and ${\mathbb{I}}_n={\left[ 0...n \right]}$, these
correspond to the classical
Ces\`aro  density and ergodic average.
We say that ${\mathbb{Y}}$ (respectively ${\mathbf{ a}}$) is {\bf stationary} if the
limit in (\ref{Cesaro.Density}) (respectively,
(\ref{Ergodic.Average})) exists and is independent of the choice of
F{\o}lner  sequence.

   Let $\ensuremath{{\mathcal{ M}}^{\sigma}\left[{\mathcal{ A}}^{\mathbb{X}}\right] }$ be the set of probability measures on
${\mathcal{ A}}^{\mathbb{X}}$ which are invariant under all ${\mathbb{X}}$-shifts,
and let $\ensuremath{{\mathcal{ M}}^{\sigma}_{e}\left[{\mathcal{ A}}^{\mathbb{X}}\right] }$ be the {\bf ergodic} measures:
\ the extremal points of 
$\ensuremath{{\mathcal{ M}}^{\sigma}\left[{\mathcal{ A}}^{\mathbb{X}}\right] }$ within $\ensuremath{{\mathcal{ M}}\left[{\mathcal{ A}}^{\mathbb{X}}\right] }$.
If $\mu\in\ensuremath{{\mathcal{ M}}^{\sigma}\left[{\mathcal{ A}}^{\mathbb{X}}\right] }$ and ${\mathbb{X}}$ is amenable, then the
generalized Birkhoff Ergodic Theorem \cite{Tempelman} says that
$\mu$-almost every ${\mathbf{ a}}\in{\mathcal{ A}}^{\mathbb{X}}$ is stationary, and, if
$\mu\in \ensuremath{{\mathcal{ M}}^{\sigma}_{e}\left[{\mathcal{ A}}^{\mathbb{X}}\right] }$, then for
$\mu$-almost all ${\mathbf{ a}}\in{\mathcal{ A}}^{\mathbb{X}}$, \ $\displaystyle{{\sf ErgAve}}_{{\mathbb{X}}}\  \phi({\mathbf{ a}}) \ = \ 
{\left\langle \phi,\mu \right\rangle }$, where we use the notational convention:
\[
  {\left\langle \phi,\mu \right\rangle } \quad = \quad \int_{{\mathcal{ A}}^{\mathbb{X}}} \phi(a_{\mathsf{ o}}) \ d\mu[{\mathbf{ a}}].
\]
  This yields the following characterization for conservation laws:

\begin{prop}
\label{Cesaro.cond}  
 Let ${\mathbb{X}}$ be an amenable group, $\phi:{\mathcal{ A}}{{\longrightarrow}}{\mathbb{R}}^+$.  The following
are equivalent:
\begin{enumerate}
  \item $\phi\in{\mathcal{ C}}({\mathfrak{ F}})$ 

  \item For any stationary ${\mathbf{ a}}\in{\mathcal{ A}}^{\mathbb{X}}$, if ${\mathbf{ a}}'={\mathfrak{ F}}({\mathbf{ a}})$, then
$\displaystyle
  \displaystyle{{\sf ErgAve}}_{{\mathbb{X}}}\  \phi({\mathbf{ a}}) =  \displaystyle{{\sf ErgAve}}_{{\mathbb{X}}}\  \phi({\mathbf{ a}}')$.

  \item For any $\mu\in\ensuremath{{\mathcal{ M}}^{\sigma}_{e}\left[{\mathcal{ A}}^{\mathbb{X}}\right] }$, if $\mu'={\mathfrak{ F}}(\mu)$, then
${\left\langle \phi,\ \mu \right\rangle } \ = \   {\left\langle \phi, \ \mu' \right\rangle }$.

  \item For any $\mu\in\ensuremath{{\mathcal{ M}}^{\sigma}\left[{\mathcal{ A}}^{\mathbb{X}}\right] }$, if $\mu'={\mathfrak{ F}}(\mu)$, then
${\left\langle \phi,\ \mu \right\rangle } \ = \   {\left\langle \phi, \ \mu' \right\rangle }$.
\end{enumerate}
 \end{prop}
\bprf
$(1\mbox{$\Longrightarrow$}2)$: Let $\{{\mathbb{I}}_n\}_{n=1}^{\infty}$ be a F{\o}lner  sequence;
if we define ${\mathbb{J}}_n = {\mathbb{I}}_n+{\mathbb{B}}$ and ${\mathbb{K}}_n =
{\mathbb{J}}_n+{\mathbb{B}}$, then $\{{\mathbb{J}}_n\}_{n=1}^{\infty}$ and
$\{{\mathbb{K}}_n\}_{n=1}^{\infty}$ are also F{\o}lner  sequences.  Let
$I_n={{\sf card}\left[{\mathbb{I}}_n\right]}$, $J_n={{\sf card}\left[{\mathbb{J}}_n\right]}$, and $K_n={{\sf card}\left[{\mathbb{K}}_n\right]}$.
Since ${\mathbb{B}}$ is
finite, the F{\o}lner  property implies: $\displaystyle \lim_{n{\rightarrow}{\infty}}
\frac{I_n}{J_n} = 1 = \lim_{n{\rightarrow}{\infty}} \frac{K_n}{J_n}$.  Given
$\epsilon>0$, find $n\in{\mathbb{N}}$ so that
\begin{eqnarray}
\label{card.bounds}
 1-\epsilon \ < \ \displaystyle\frac{I_n}{J_n} &\mbox{and}&
\displaystyle\frac{K_n}{J_n} \ <\  1+\epsilon; \\
\label{low.bound}
 (1-\epsilon)\cdot \displaystyle{{\sf ErgAve}}_{{\mathbb{X}}}\ \phi({\mathbf{ a}})
& <&\displaystyle\frac{1}{I_n}\sum_{\ensuremath{\mathsf{ i}}\in{\mathbb{I}}_n} \phi(a_\ensuremath{\mathsf{ i}}); \\ 
\label{high.bound}
\mbox{and} \ \ \ \ 
\displaystyle\frac{1}{K_n}\sum_{{\mathsf{ k}}\in{\mathbb{K}}_n} \phi(a_{\mathsf{ k}})
& <& (1+\epsilon)\cdot  \displaystyle{{\sf ErgAve}}_{{\mathbb{X}}}\ \phi({\mathbf{ a}}).
\end{eqnarray}
Now, ${\mathbb{I}}_n\subset{{\mathbf{ i}}{\mathbf{ n}}{\mathbf{ t}}\left[ {\mathbb{J}}_n \right]}$ and ${\mathbb{K}}_n={{\mathbf{ c}}{\mathbf{ l}}\left[ {\mathbb{J}}_n \right]}$, so applying
Theorem \ref{inf.cond} to ${\mathbb{J}}_n$ yields:
\begin{equation}
\label{sum.bounds}
\sum_{\ensuremath{\mathsf{ i}}\in{\mathbb{I}}_n} \phi(a_\ensuremath{\mathsf{ i}}) 
\ \ \leq \ \
\sum_{{\mathsf{ j}}\in{\mathbb{J}}_n} \phi(a'_{\mathsf{ j}})
\ \ \leq \ \
\sum_{{\mathsf{ k}}\in{\mathbb{K}}_n} \phi(a_{\mathsf{ k}}).
\end{equation}
\begin{eqnarray*}
\mbox{Thus,} \ \ 
(1-\epsilon)^2 \cdot \displaystyle{{\sf ErgAve}}_{{\mathbb{X}}}\ \phi({\mathbf{ a}})
&<_{[\ref{low.bound}]}&
\frac{1-\epsilon}{I_n} \sum_{\ensuremath{\mathsf{ i}}\in{\mathbb{I}}_n} \phi(a_\ensuremath{\mathsf{ i}}) 
\ <_{[\ref{card.bounds}]} \ 
\frac{1}{J_n} \sum_{\ensuremath{\mathsf{ i}}\in{\mathbb{I}}_n} \phi(a_\ensuremath{\mathsf{ i}}) \\
& \leq_{[\ref{sum.bounds}]} &
\frac{1}{J_n} \sum_{{\mathsf{ j}}\in{\mathbb{J}}_n} \phi(a'_{\mathsf{ j}}) \\
& \leq_{[\ref{sum.bounds}]} &
\frac{1}{J_n} \sum_{{\mathsf{ k}}\in{\mathbb{K}}_n} \phi(a_{\mathsf{ k}}) 
\ <_{[\ref{card.bounds}]} \
\frac{1+\epsilon}{K_n} \sum_{{\mathsf{ k}}\in{\mathbb{K}}_n} \phi(a_{\mathsf{ k}}) \\
& <_{[\ref{high.bound}]} & 
(1+\epsilon)^2 \cdot \displaystyle{{\sf ErgAve}}_{{\mathbb{X}}}\ \phi({\mathbf{ a}}),
\end{eqnarray*}
 where each inequality follows from formula with the same number.
Letting $\epsilon{\rightarrow} 0$ as $n{\rightarrow}{\infty}$, we conclude by a squeezing
argument:
\[
 \displaystyle{{\sf ErgAve}}_{{\mathbb{X}}}\ \phi({\mathbf{ a}}) 
\quad = \quad  \lim_{\epsilon\rightarrow 0} (1\pm\epsilon)^2 \cdot  \displaystyle{{\sf ErgAve}}_{{\mathbb{X}}}\ \phi({\mathbf{ a}})
\quad = \quad 
  \lim_{n{\rightarrow}{\infty}} \frac{1}{J_n} \sum_{{\mathsf{ j}}\in{\mathbb{J}}_n} \phi(a'_{\mathsf{ j}}) 
\quad = \quad  \displaystyle{{\sf ErgAve}}_{{\mathbb{X}}}\ \phi({\mathbf{ a}}').
\]
$(2\mbox{$\Longrightarrow$}1)$:  Let ${\mathbf{ b}}\in{\mathcal{ A}}^{<{\mathbb{X}}}$, with ${\mathbf{ b}}'={\mathfrak{ F}}({\mathbf{ b}})$.
Let ${\mathbb{U}}={\sf supp}\left[{\mathbf{ b}}\right]$ and  ${\mathbb{U}}'={{\mathbf{ c}}{\mathbf{ l}}\left[ {\mathbb{U}} \right]}$,
and find stationary ${\mathbb{Y}}\subset{\mathbb{X}}$ with ${{\sf density}\left[{\mathbb{Y}}\right] }=\delta>0$, 
such that ${\mathsf{ y}}_1+{\mathbb{U}}'$ and
${\mathsf{ y}}_2+{\mathbb{U}}'$ are disjoint for any ${\mathsf{ y}}_1\neq{\mathsf{ y}}_2\in{\mathbb{Y}}$.  Then
define ${\mathbf{ a}}\in{\mathcal{ A}}^{\mathbb{X}}$ by: ${\mathbf{ a}}\raisebox{-0.3em}{$\left|_{{\mathsf{ y}}+{\mathbb{U}}}\right.$}= {{{\boldsymbol{\sigma}}}^{{\mathsf{ y}}}} ({\mathbf{ b}})$
for every ${\mathsf{ y}}\in{\mathbb{Y}}$, and $a_{\mathsf{ x}}=\mbox{\cursive O}\in\mathbf{0}$ for all ${\mathsf{ x}}\not\in{\mathbb{Y}}+{\mathbb{U}}$.
Let ${\mathbf{ a}}'={\mathfrak{ F}}({\mathbf{ a}})$.

\refstepcounter{claimcount}                {\bf Claim \theclaimcount: \ }{\sl  $ \displaystyle{{\sf ErgAve}}_{{\mathbb{X}}}\  \phi({\mathbf{ a}}) \ = \ \delta\cdot \Sigma\phi({\mathbf{ b}})$
and 
$ \displaystyle{{\sf ErgAve}}_{{\mathbb{X}}}\  \phi({\mathbf{ a}}') \ =  \ \delta\cdot \Sigma\phi({\mathbf{ b}}')$.

}
\bprf 
Let $\{{\mathbb{I}}_n\}_{n=1}^{\infty}$ be a F{\o}lner  sequence; \  for any 
$n\in{\mathbb{N}}$, let ${\mathbb{Y}}_n = {\mathbb{Y}}\cap{\mathbb{I}}_n$ and ${\mathbb{Y}}^\ast_n = 
{\left\{ {\mathsf{ y}}\in{\mathbb{Y}} \; ; \; {\mathsf{ y}}+{\mathbb{U}}\subset{\mathbb{I}}_n \right\} }$. 
Assume ${\mathbb{U}}$ contains ${\mathsf{ O}}$, the identity element of ${\mathbb{X}}$; \ thus 
${\mathbb{Y}}^\ast_n \subset {\mathbb{Y}}_n$. 
By construction, for any ${\mathsf{ y}}\in{\mathbb{Y}}$, \ \
$\displaystyle  \sum_{{\mathsf{ u}}\in{\mathbb{U}}} \phi\left(a_{{\mathsf{ y}}+{\mathsf{ u}}}\right) = \Sigma\phi({\mathbf{ b}})$. \ \
Thus,
\begin{eqnarray*}
{{\sf card}\left[{\mathbb{Y}}^\ast_n\right]} \cdot \Sigma\phi({\mathbf{ b}}) 
& = & 
 \sum_{{\mathsf{ y}}\in{\mathbb{Y}}^\ast_n} \ \sum_{{\mathsf{ u}}\in{\mathbb{U}}} \phi\left(a_{{\mathsf{ y}}+{\mathsf{ u}}}\right) \\
& \leq &
 \sum_{\ensuremath{\mathsf{ i}}\in{\mathbb{I}}_n} \phi(a_\ensuremath{\mathsf{ i}})
\ \ \leq \ \
 \sum_{{\mathsf{ y}}\in{\mathbb{Y}}_n} \ \sum_{{\mathsf{ u}}\in{\mathbb{U}}} \phi\left(a_{{\mathsf{ y}}+{\mathsf{ u}}}\right)
\ \ = \ \
{{\sf card}\left[{\mathbb{Y}}_n\right]} \cdot \Sigma\phi({\mathbf{ b}})
\end{eqnarray*}
Now divide everything by $I_n={{\sf card}\left[{\mathbb{I}}_n\right]}$, and take the limit as $n{\rightarrow}{\infty}$.
By definition, ${{\sf density}\left[{\mathbb{Y}}\right] }=\delta$, so that
\[
 \lim_{n{\rightarrow}{\infty}} \frac{{{\sf card}\left[{\mathbb{Y}}^\ast_n\right]}}{I_n} \quad = \quad 
\delta \quad  = \quad 
\lim_{n{\rightarrow}{\infty}} \frac{{{\sf card}\left[{\mathbb{Y}}_n\right]}}{I_n};
\]
  thus, by a squeezing argument, 
\[
\displaystyle{{\sf ErgAve}}_{{\mathbb{X}}}\  \phi({\mathbf{ a}})
\quad = \quad
\lim_{n{\rightarrow}{\infty}} \frac{1}{I_n} \sum_{\ensuremath{\mathsf{ i}}\in{\mathbb{I}}_n} \phi(a_\ensuremath{\mathsf{ i}})
\quad = \quad
  \delta\cdot \Sigma\phi({\mathbf{ b}}).
\] 
 The proof for ${\mathbf{ a}}'$ and ${\mathbf{ b}}'$  uses ${\mathbb{U}}'$ instead
 of ${\mathbb{U}}$, and the fact that, for any ${\mathsf{ y}}\in{\mathbb{Y}}$, \ \ 
$\displaystyle \sum_{{\mathsf{ u}}'\in{\mathbb{U}}'} \phi\left(a'_{{\mathsf{ y}}+{\mathsf{ u}}'}\right) = \Sigma\phi({\mathbf{ b}}')$.
 {\tt \dotfill~$\Box$~[Claim~\theclaimcount] }\end{list} 
Thus,
$
\delta\cdot \Sigma\phi({\mathbf{ b}})
\  =  \ 
  \displaystyle{{\sf ErgAve}}_{{\mathbb{X}}}\  \phi({\mathbf{ a}}) 
\  =_{[2]}  \ 
  \displaystyle{{\sf ErgAve}}_{{\mathbb{X}}}\  \phi({\mathbf{ a}}') 
\  =  \ 
 \delta\cdot \Sigma\phi({\mathbf{ b}}')$, \
 where inequality $[2]$ follows from hypothesis (2).  This implies
$\Sigma\phi({\mathbf{ b}}) = \Sigma\phi({\mathbf{ b}}')$; \ since this holds for any ${\mathbf{ b}}\in{\mathcal{ A}}^{<{\mathbb{X}}}$, we conclude that $\phi\in{\mathcal{ C}}({\mathfrak{ F}})$.

\medskip

{$(2\mbox{$\Longrightarrow$}3)$:} Apply the Birkhoff Ergodic Theorem.

{$(3\mbox{$\Longrightarrow$}4)$:} Any element of $\ensuremath{{\mathcal{ M}}^{\sigma}\left[{\mathcal{ A}}^{\mathbb{X}}\right] }$,
 is a weak*-limit of convex combinations of ergodic measures. 
So, suppose $\displaystyle\mu=\mathrm{wk^*}\!\!\lim_{i{\rightarrow}{\infty}}\nu_i$, where $\displaystyle\nu_i =
\sum_{j=1}^{J_i} \lambda_{i j} \, \eta_{i j}$, with $\eta_{i
j}\in\ensuremath{{\mathcal{ M}}^{\sigma}_{e}\left[{\mathcal{ A}}^{\mathbb{X}}\right] }$, and $\lambda_{i j}\in{\left[ 0,1 \right]}$, for all
$i\in{\mathbb{N}}$ and $j\in{\left[ 1...J_i \right]}$.  If $\mu'={\mathfrak{ F}}(\mu)$, then
$\displaystyle\mu'=\mathrm{wk^*}\!\!\lim_{i{\rightarrow}{\infty}}\nu'_i$, where $\displaystyle \nu'_i={\mathfrak{ F}}(\nu_i) =
\sum_{j=1}^{J_i} \lambda_{i j}\,\eta'_{i j}$, with $\eta_{i j}' =
{\mathfrak{ F}}(\eta_{i j})$.
But by hypothesis (3), ${\left\langle \phi,\eta_{i j} \right\rangle } = {\left\langle \phi,\eta_{i
j}' \right\rangle }$ for all $i$ and $j$.  Thus, ${\left\langle \phi,\nu_{i} \right\rangle } =
{\left\langle \phi,\nu_{i}' \right\rangle }$ for all $i\in{\mathbb{N}}$; thus ${\left\langle \phi,\mu \right\rangle } =
{\left\langle \phi,\mu' \right\rangle }$.

\paragraph*{$(4\mbox{$\Longrightarrow$}2)$:}  If ${\mathbf{ a}}\in{\mathcal{ A}}^{\mathbb{X}}$ is stationary,
let $\delta_{\mathbf{ a}}\in\ensuremath{{\mathcal{ M}}\left[{\mathcal{ A}}^{\mathbb{X}}\right] }$ be the point mass at ${\mathbf{ a}}$;
then $\delta_{{\mathbf{ a}}'} = {\mathfrak{ F}}(\delta_{\mathbf{ a}})$.  Let $\{{\mathbb{I}}_n\}_{n=1}^{\infty}$ be
a F{\o}lner  sequence, and for all $n\in{\mathbb{N}}$,
let $\displaystyle\mu_n = \frac{1}{I_n} \sum_{\ensuremath{\mathsf{ i}}\in{\mathbb{I}}_n}  {{{\boldsymbol{\sigma}}}^{\ensuremath{\mathsf{ i}}}} \delta_{\mathbf{ a}}$ and
$\displaystyle\mu'_n = \frac{1}{I_n} \sum_{\ensuremath{\mathsf{ i}}\in{\mathbb{I}}_n}  {{{\boldsymbol{\sigma}}}^{\ensuremath{\mathsf{ i}}}} \delta_{{\mathbf{ a}}'}$.
Since $\ensuremath{{\mathcal{ M}}\left[{\mathcal{ A}}^{\mathbb{X}}\right] }$ is compact in the weak* topology, the
sequence $\{\mu_n\}_{n=1}^{\infty}$ has a weak* cluster point, $\mu$,
which  by construction is shift-invariant.
Dropping to a subsequence if necessary, we'll say
$\displaystyle\mu = \mathrm{wk^*}\!\!\lim_{n{\rightarrow}{\infty}} \mu_n$.  Thus,
\[
{\left\langle \phi,\mu \right\rangle }
\ = \ 
\lim_{n{\rightarrow}{\infty}}
 {\left\langle \phi,\mu_n \right\rangle } \ = \ 
\lim_{n{\rightarrow}{\infty}} \frac{1}{I_n} \sum_{\ensuremath{\mathsf{ i}}\in{\mathbb{I}}_n} 
 {\left\langle \phi, \  {{{\boldsymbol{\sigma}}}^{\ensuremath{\mathsf{ i}}}} \delta_{\mathbf{ a}} \right\rangle }
\ = \ 
\lim_{n{\rightarrow}{\infty}} \frac{1}{I_n} \sum_{\ensuremath{\mathsf{ i}}\in{\mathbb{I}}_n} \phi(a_\ensuremath{\mathsf{ i}}) 
\ = \ \displaystyle{{\sf ErgAve}}_{{\mathbb{X}}}\  \phi({\mathbf{ a}}).
\]
If $\mu'={\mathfrak{ F}}(\mu)$, then $\displaystyle\mu' = \mathrm{wk^*}\!\!\lim_{n{\rightarrow}{\infty}} \mu'_n$
is also shift-invariant, and ${\left\langle \phi,\mu' \right\rangle } \ = \  \displaystyle{{\sf ErgAve}}_{{\mathbb{X}}}\  \phi({\mathbf{ a}}')$.  But by hypothesis (4), we have ${\left\langle \phi,\mu' \right\rangle } = {\left\langle \phi,\mu \right\rangle }$;
hence,  $\displaystyle{{\sf ErgAve}}_{{\mathbb{X}}}\  \phi({\mathbf{ a}}') =  \displaystyle{{\sf ErgAve}}_{{\mathbb{X}}}\  \phi({\mathbf{ a}})$.
 {\tt \hrulefill $\Box$ } \end{list}  \medskip  

  Real-valued conservation laws thus preclude unique ergodicity:

\begin{cor}  
 Let $\displaystyle R_0=\min_{a\in{\mathcal{ A}}} \phi(a)$ and $\displaystyle R_1=\max_{a\in{\mathcal{ A}}} \phi(a)$.
If $\phi\in{\mathcal{ C}}({\mathfrak{ F}})$,
then for any $r\in{\left[ R_0,R_1 \right]}$, there is an ${\mathfrak{ F}}$-invariant measure
$\mu_r\in\ensuremath{{\mathcal{ M}}^{\sigma}\left[{\mathcal{ A}}^{\mathbb{X}}\right] }$ such that
${\left\langle \phi, \ \mu_r \right\rangle } \ = \  r$.
 \end{cor}
\bprf
  Let $a_k\in{\mathcal{ A}}$ be such that $\phi(a_k)=R_k$.  Given
$r\in{\left[ R_0,R_1 \right]}$, let $\lambda\in{\left[ 0,1 \right]}$ be such that $r=\lambda R_0 +
(1-\lambda) R_1$.  Let $\rho$ be the probability measure on ${\mathcal{ A}}$ with
$\rho\{a_0\}=\lambda$ and $\rho\{a_1\}=1-\lambda$, and let $\nu_r$ be the
associated Bernoulli measure on ${\mathcal{ A}}^{\mathbb{X}}$ ---that is, the product
measure $\nu_r = \displaystyle\bigotimes_{{\mathsf{ x}}\in{\mathbb{X}}} \rho$.

Thus $\nu_r$ is shift-ergodic and ${\left\langle \phi, \ \nu_r \right\rangle } \ = \
r$.
  For all $N\in{\mathbb{N}}$, define $\displaystyle \eta_{N} = \frac{1}{N}\sum_{n=1}^N
{\mathfrak{ F}}^n \nu_r$; then by part 3 of Proposition \ref{Cesaro.cond}, ${\left\langle \phi, \ \eta_{N} \right\rangle } \ = \
r$.  Since $\ensuremath{{\mathcal{ M}}^{\sigma}\left[{\mathcal{ A}}^{\mathbb{X}}\right] }$ is compact in the weak* topology, the
sequence $\{\eta_N\}_{n=1}^{\infty}$ has a weak* limit point, $\mu_r$.  By
construction, $\mu_r$ is ${\mathfrak{ F}}$-invariant, shift-invariant, and 
${\left\langle \phi, \ \mu_r \right\rangle } \ = \  r$. 
 {\tt \hrulefill $\Box$ } \end{list}  \medskip

\section{Characterizations by Measure
\label{S:measure}}

  If ${\mathbb{X}}$ is not amenable, then the methods of \S\ref{S:Cesaro}
are inapplicable.  However, we can still characterize conservation laws
on ${\mathcal{ A}}^{\mathbb{X}}$ in terms of shift-invariant
measures, by projecting ${\mathbb{X}}$ onto a finite quotient group.

   Let ${\mathfrak{ Q}}$ be as in \S\ref{S:finitary}.  If ${\widetilde{\mathbb{X}}}\in{\mathfrak{ Q}}$, then 
let $\ensuremath{{\mathcal{ M}}^{\sigma}\left[{\mathcal{ A}}^{\widetilde{\mathbb{X}}};\ {\mathbb{R}}\right] }$ be the space of shift-invariant,
real-valued measures on ${\mathcal{ A}}^{\widetilde{\mathbb{X}}}$, and $\ensuremath{{\mathcal{ M}}^{\sigma}\left[{\mathcal{ A}}^{\widetilde{\mathbb{X}}}\right] }\subset
\ensuremath{{\mathcal{ M}}^{\sigma}\left[{\mathcal{ A}}^{\widetilde{\mathbb{X}}}; \ {\mathbb{R}}\right] }$ the space of shift-invariant probability
measures.   Let ${\widetilde{\mathbb{B}}}\subset{\widetilde{\mathbb{X}}}$ be the (bijective) image of ${\mathbb{B}}$
under the quotient map ${\mathbb{X}}{{\longrightarrow}}{\widetilde{\mathbb{X}}}$, and,
for any $\mu\in\ensuremath{{\mathcal{ M}}^{\sigma}\left[{\mathcal{ A}}^{\widetilde{\mathbb{X}}};\ {\mathbb{R}}\right] }$, let $\mu_{\widetilde{\mathbb{B}}} =
{\mathbf{ pr}_{{{\widetilde{\mathbb{B}}}}}}^* (\mu) \in \ensuremath{{\mathcal{ M}}\left[{\mathcal{ A}}^{\widetilde{\mathbb{B}}};\ {\mathbb{R}}\right] }$ be the marginal
projection of $\mu$ onto ${\widetilde{\mathbb{B}}}$.  Via the bijection ${\mathbb{B}}{{\longrightarrow}}{\widetilde{\mathbb{B}}}$,
we can identify ${\mathcal{ A}}^{\widetilde{\mathbb{B}}}$ with ${\mathcal{ A}}^{\mathbb{B}}$, and
thus $\ensuremath{{\mathcal{ M}}\left[{\mathcal{ A}}^{\widetilde{\mathbb{B}}};\ {\mathbb{R}}\right] }$ with 
$\ensuremath{{\mathcal{ M}}\left[{\mathcal{ A}}^{\mathbb{B}};\ {\mathbb{R}}\right] }$.  Now define:
\[
\ensuremath{{\mathcal{ M}}^{\sigma}\left[{\mathcal{ A}}^{\mathbb{B}};\ {\mathbb{R}}\right] }
 \ =\  {\left\{ \nu\in\ensuremath{{\mathcal{ M}}\left[{\mathcal{ A}}^{\mathbb{B}};{\mathbb{R}}\right] } \; ; \; \nu=\mu_{\widetilde{\mathbb{B}}}
\ \mbox{for some} \ \mu\in\ensuremath{{\mathcal{ M}}^{\sigma}\left[{\mathcal{ A}}^{\widetilde{\mathbb{X}}}; \ {\mathbb{R}}\right] } \ \mbox{and }
\ {\widetilde{\mathbb{X}}}\in{\mathfrak{ Q}} \right\} }.
\]
  Define ${\eth{\phi}}:{\mathcal{ A}}^{\mathbb{B}}{{\longrightarrow}}{\mathbb{R}}$
by \ ${\eth{\phi}\left({\mathbf{ a}}\right)} = \phi\left({\mathfrak{ f}}({\mathbf{ a}})\right) - \phi(a_{\mathsf{ o}})$, for
any ${\mathbf{ a}}\in{\mathcal{ A}}^{\mathbb{B}}$.  For any ${\widetilde{\mathbb{X}}}\in{\mathfrak{ Q}}$, define
${\partial_t\, {\phi}_{}}:{\mathcal{ A}}^{\widetilde{\mathbb{X}}}{{\longrightarrow}}{\mathbb{R}}^{\widetilde{\mathbb{X}}}$ by
${\partial_t\, {\phi}_{{\mathsf{ x}}}\left({\mathbf{ a}}\right)} = \phi\left({\mathfrak{ F}}({\mathbf{ a}})_{\mathsf{ x}}\right) - \phi(a_{\mathsf{ x}})$, \ 
for all ${\mathbf{ a}}\in{\mathcal{ A}}^{\widetilde{\mathbb{X}}}$ and ${\mathsf{ x}}\in{\widetilde{\mathbb{X}}}$.

\begin{prop}\label{measure.thm}   Let $\phi:{\mathcal{ A}}{{\longrightarrow}}{\mathbb{R}}$ and let ${\mathfrak{ F}}$ be a CA.
 The following are equivalent:
  
\begin{enumerate}
 \item  \label{M1} $\phi\in{\mathcal{ C}}\left({\mathfrak{ F}}\right)$.

 \item  \label{M2} For all $ {\widetilde{\mathbb{X}}}\in{\mathfrak{ Q}}$, \  $\phi\in{\mathcal{ C}}\left({\widetilde{\mathfrak{ F}}}\right)$.

 \item  \label{M3}  For all $ {\widetilde{\mathbb{X}}}\in{\mathfrak{ Q}}$, and all $ 
\mu\in\ensuremath{{\mathcal{ M}}^{\sigma}\left[{\mathcal{ A}}^{\widetilde{\mathbb{X}}}; \ {\mathbb{R}}\right] }$, \ ${\left\langle \phi_{\mathsf{ o}}, \ \mu \right\rangle }=
 {\left\langle \phi_{\mathsf{ o}}, \ {\mathfrak{ F}}(\mu) \right\rangle }$.

 \item  \label{M4} For all $ {\widetilde{\mathbb{X}}}\in{\mathfrak{ Q}}$, and all $ 
\mu\in\ensuremath{{\mathcal{ M}}^{\sigma}\left[{\mathcal{ A}}^{\widetilde{\mathbb{X}}}; \ {\mathbb{R}}\right] }$, \ ${\left\langle {\partial_t\, {\phi}_{{\mathsf{ o}}}}, \ \mu \right\rangle }=0$.

 \item  \label{M5} For all $ \mu_{\widetilde{\mathbb{B}}}\in\ensuremath{{\mathcal{ M}}^{\sigma}\left[{\mathcal{ A}}^{\mathbb{B}}; \ {\mathbb{R}}\right] }$, \
 ${\left\langle {\eth{\phi}}, \ \mu_{\widetilde{\mathbb{B}}} \right\rangle }=0$.
\end{enumerate}
 \end{prop}
\bprf  $(\ref{M1})\Leftrightarrow (\ref{M2})$:  This just restates
Corollary \ref{fuks.boccara}.

 $(\ref{M4}) \Leftrightarrow (\ref{M5})$:  If $\mu\in \ensuremath{{\mathcal{ M}}^{\sigma}\left[{\mathcal{ A}}^{\widetilde{\mathbb{X}}}; \ {\mathbb{R}}\right] }$, then  ${\left\langle {\partial_t\, {\phi}_{{\mathsf{ o}}}}, \ \mu \right\rangle } = {\left\langle {\eth{\phi}},\ \mu_{\widetilde{\mathbb{B}}} \right\rangle }$.

$(\ref{M3}) \Leftrightarrow (\ref{M4})$: By definition, 
${\left\langle {\partial_t\, {\phi}_{{\mathsf{ o}}}},\ \mu \right\rangle } =
{\left\langle \phi_{\mathsf{ o}}\circ{\mathfrak{ F}}-\phi_{\mathsf{ o}}, \ \mu \right\rangle } =
{\left\langle \phi_{\mathsf{ o}},\ {\mathfrak{ F}}(\mu) \right\rangle } - {\left\langle \phi_{\mathsf{ o}},\ \mu \right\rangle }$.

$(\ref{M2}) \mbox{$\Longrightarrow$} (\ref{M4})$:\ ${\widetilde{\mathbb{X}}}$ is finite, so
${\mathcal{ A}}^{<{\widetilde{\mathbb{X}}}} = {\mathcal{ A}}^{\widetilde{\mathbb{X}}}$. 
 Thus, we can well-define $\overline{{\partial_t\, {\phi}_{}}}:{\mathcal{ A}}^{\widetilde{\mathbb{X}}}{{\longrightarrow}}{\mathbb{R}}$ by:
$\displaystyle\overline{{\partial_t\, {\phi}_{}}}({\mathbf{ a}})=\sum_{{\mathsf{ x}}\in{\widetilde{\mathbb{X}}}} {\partial_t\, {\phi}_{{\mathsf{ x}}}\left({\mathbf{ a}}\right)}$.
If $\mu\in\ensuremath{{\mathcal{ M}}\left[{\mathcal{ A}}^{\widetilde{\mathbb{X}}}; \ {\mathbb{R}}\right] }$, and we likewise define $\displaystyle\overline{\mu}=
\sum_{{\mathsf{ x}}\in{\widetilde{\mathbb{X}}}} {{{\boldsymbol{\sigma}}}^{{\mathsf{ x}}}} \mu$, then clearly, \ 
$\displaystyle {\left\langle {\partial_t\, {\phi}_{{\mathsf{ o}}}}, \ \overline{\mu} \right\rangle } \ = \
{\left\langle \overline{{\partial_t\, {\phi}_{}}},\mu \right\rangle }$.
 
But if $\phi\in{\mathcal{ C}}\left({\widetilde{\mathfrak{ F}}}\right)$, then $\overline{{\partial_t\, {\phi}_{}}}\equiv0$.
To see this, let ${\mathbf{ a}}\in{\mathcal{ A}}^{\widetilde{\mathbb{X}}}$; then:
\[
\overline{{\partial_t\, {\phi}_{}}}({\mathbf{ a}})
\quad=\quad  
\sum_{{\mathsf{ x}}\in{\widetilde{\mathbb{X}}}} \left(\rule[-0.5em]{0em}{1em}
 \phi_{\mathsf{ x}}\,{\mathfrak{ F}}({\mathbf{ a}})  - \phi_{\mathsf{ x}}({\mathbf{ a}})\right)
\quad=\quad  
\Sigma\phi\,{\mathfrak{ F}}({\mathbf{ a}})
-
\Sigma\phi({\mathbf{ a}}) 
\quad  = \quad 
0.
\]
  Also, if $\mu\in\ensuremath{{\mathcal{ M}}^{\sigma}\left[{\mathcal{ A}}^{\widetilde{\mathbb{X}}};\ {\mathbb{R}}\right] }$, then $\overline{\mu}={{\sf card}\left[{\widetilde{\mathbb{X}}}\right]}\cdot\mu$.
Combining these facts yields:
\[
{{\sf card}\left[{\widetilde{\mathbb{X}}}\right]}\cdot {\left\langle {\partial_t\, {\phi}_{{\mathsf{ o}}}}, \ \mu \right\rangle }
\quad=\quad  
{\left\langle {\partial_t\, {\phi}_{{\mathsf{ o}}}}, \ \overline{\mu} \right\rangle }
\quad=\quad 
{\left\langle \overline{{\partial_t\, {\phi}_{}}},\mu \right\rangle } 
\quad=\quad  
{\left\langle 0,\mu \right\rangle } \quad=\quad  0. 
\] 

{$(\ref{M4})\mbox{$\Longrightarrow$} (\ref{M2})$}:   If ${\mathbf{ a}}\in{\mathcal{ A}}^{\widetilde{\mathbb{X}}}$,
and ${\boldsymbol{\delta}}_{\mathbf{ a}}\in\ensuremath{{\mathcal{ M}}\left[{\mathcal{ A}}^{\widetilde{\mathbb{X}}}\right] }$ is the point mass at ${\mathbf{ a}}$, then $\displaystyle
\overline{{\boldsymbol{\delta}}_{\mathbf{ a}}} = \sum_{{\mathsf{ x}}\in{\widetilde{\mathbb{X}}}}
 {{{\boldsymbol{\sigma}}}^{{\mathsf{ x}}}} {\boldsymbol{\delta}}_{\mathbf{ a}}$ is in $\ensuremath{{\mathcal{ M}}^{\sigma}\left[{\mathcal{ A}}^{\widetilde{\mathbb{X}}}; \ {\mathbb{R}}\right] }$.  Thus,
\begin{eqnarray*}
 \Sigma\phi\,{\mathfrak{ F}}({\mathbf{ a}}) - \Sigma\phi({\mathbf{ a}})
& = &
 \sum_{{\mathsf{ x}}\in{\widetilde{\mathbb{X}}}}
      \left(\rule[-0.5em]{0em}{1em} \phi_{\mathsf{ x}}\,{\mathfrak{ F}}({\mathbf{ a}}) - \phi_{\mathsf{ x}}({\mathbf{ a}})\right)
\ \quad=\quad  \
\overline{{\partial_t\, {\phi}_{}}}({\mathbf{ a}})\\
& = &
{\left\langle \overline{{\partial_t\, {\phi}_{}}}, \ {\boldsymbol{\delta}}_{\mathbf{ a}} \right\rangle }
\ \quad=\quad  \
{\left\langle {\partial_t\, {\phi}_{}}, \ \overline{{\boldsymbol{\delta}}_{\mathbf{ a}}} \right\rangle }
\ \quad=\quad  \ 0.
\end{eqnarray*}
This is true for any ${\mathbf{ a}}\in{\mathcal{ A}}^{\widetilde{\mathbb{X}}}$, so
$\phi\in{\mathcal{ C}}\left({\widetilde{\mathfrak{ F}}}\right)$.
 {\tt \hrulefill $\Box$ } \end{list}  \medskip  

  ${\mathcal{ A}}^{\mathbb{B}}$ is finite, so $\ensuremath{{\mathcal{ M}}\left[{\mathcal{ A}}^{\mathbb{B}}; \ {\mathbb{R}}\right] }$ is a finite
dimensional vector space, and $\ensuremath{{\mathcal{ M}}^{\sigma}\left[{\mathcal{ A}}^{\mathbb{B}}; \ {\mathbb{R}}\right] }$ is a
linear subspace, with some finite basis $\mu_1,\ldots,\mu_N$
(for example, see \cite{PivatoCJM}). To check
{\bf Part 5} of Proposition \ref{measure.thm}, it suffices to check
that ${\left\langle {\eth{\phi}},\ \mu_n \right\rangle } = 0$ for all $n\in{\left[ 1..N \right]}$,
a finite system of linear equations.  

\begin{cor}\label{measure.cor}  Let ${{\sf card}\left[{\mathcal{ A}}\right]}=A$ and ${{\sf card}\left[{\mathbb{B}}\right]}=B$.  If
$\phi\in{\mathcal{ C}}({\mathfrak{ F}})$, then $\displaystyle\sum_{{\mathbf{ a}}\in{\mathcal{ A}}^{\mathbb{B}}}
\phi\circ{\mathfrak{ f}}({\mathbf{ a}}) \ = \ A^{B-1} \sum_{a\in{\mathcal{ A}}} \phi(a)$. \ \ 
In particular, if ${\mathcal{ A}}={\left[ 0...A \right)}\subset{\mathbb{N}}$,
 and $\phi(a)=a$ for all $a\in{\mathcal{ A}}$,
then \  $\displaystyle\sum_{{\mathbf{ a}}\in{\mathcal{ A}}^{\mathbb{B}}} {\mathfrak{ f}}({\mathbf{ a}}) \ = \ \frac{1}{2} A^B (A-1)$. \end{cor}
\bprf
Let $\eta\in\ensuremath{{\mathcal{ M}}^{\sigma}\left[{\mathcal{ A}}^{\mathbb{B}}\right] }$ be the uniform Bernoulli measure,
assigning probability $\displaystyle\frac{1}{A^B}$ to
every element of ${\mathcal{ A}}^{\mathbb{B}}$.  Then
$\displaystyle\frac{1}{A^B}\sum_{{\mathbf{ a}}\in{\mathcal{ A}}^{\mathbb{B}}} \phi\circ{\mathfrak{ f}}({\mathbf{ a}})
- \frac{1}{A}\sum_{a\in{\mathcal{ A}}} \phi(a)
\quad = \quad  {\left\langle \phi\circ{\mathfrak{ f}}, \ \eta \right\rangle }- {\left\langle \phi, \ \eta \right\rangle }
\quad = \quad  {\left\langle {\eth{\phi}}, \ \eta \right\rangle }
\quad = \quad 0$,
by {\bf Part 5} of  Proposition \ref{measure.thm}.   Thus,
$\displaystyle\frac{1}{A^B}\sum_{{\mathbf{ a}}\in{\mathcal{ A}}^{\mathbb{B}}} \phi\circ{\mathfrak{ f}}({\mathbf{ a}}) \ = \ 
 \frac{1}{A}\sum_{a\in{\mathcal{ A}}} \phi(a).$
 {\tt \hrulefill $\Box$ } \end{list}  \medskip  

 The second statement of Corollary \ref{measure.cor} generalizes
Corollary 2.1 of \cite{BoccaraFuks2}, which is the analogous
result when ${\mathbb{X}}={\mathbb{Z}}$ and ${\mathbb{B}}$ is an interval.
\section{Constructing Conservative Cellular Automata
\label{S:construct}}

  Given $\phi:{\mathcal{ A}}{{\longrightarrow}}{\mathbb{R}}$, can we construct a cellular automaton
${\mathfrak{ F}}:{\mathcal{ A}}^{\mathbb{X}}\!\longrightarrow{\mathcal{ A}}^{\mathbb{X}}\!$ \ so that $\phi\in{\mathcal{ C}}({\mathfrak{ F}})$?  Can we enumerate
all such automata?  By {\bf Part 2} of Proposition \ref{recode.thm},
we can assume $\phi:{\mathcal{ A}}{{\longrightarrow}}{\mathbb{N}}^K$; hence, the problem is to
construct a CA which preserves $K$ species of particles.
For simplicity, we will consider the construction of a CA
preserving one species of particle.

  PPCA are usually constructed by explicitly specifying how each
particle displaces from its current position to a nearby location, in
response to its current local environment.  The local rule of the
PPCA, as a map ${\mathfrak{ f}}:{\mathcal{ A}}^{\mathbb{B}}{{\longrightarrow}}{\mathcal{ A}}$, is then formulated {\em a
posteori} to realize this ``displacement'' model.  Does every PPCA
arise in this manner?  Or are there PPCA not admitting any
displacement representation?

  In the case ${\mathbb{X}}={\mathbb{Z}}$, ${\mathcal{ A}}={\left[ 0...A \right]}$, and $\phi(a)=a$, every
PPCA unearthed in extensive computational searches has admitted a
displacement representation \cite{BoccaraFuks1,BoccaraFuks2}.  But it
is not clear why such a representation should always exist, and, for
more complicated  PPCA, such a representation, even if it exists, may not
be obvious from inspection.

  If a displacement rule is to yield a cellular automaton, it must
satisfy the following conditions:

\begin{description}
\item[(D1)] The  rule is {\em equivariant}  under shifts: \
a particle at ${\mathsf{ x}}\in{\mathbb{X}}$, in configuration
${\mathbf{ a}}\in{\mathcal{ A}}^{\mathbb{X}}$ will experience the same displacement as a
particle at ${\mathsf{ x}}+{\mathsf{ y}}$ in configuration $ {{{\boldsymbol{\sigma}}}^{{\mathsf{ y}}}} ({\mathbf{ a}})$.

\item[(D2)] Each particle has {\em bounded velocity}.
The new position of any particle at ${\mathsf{ x}}$ is inside ${\mathsf{ x}}+{{\mathbb{B}}}$.

\item[(D3)]  Each particle's displacement {\em locally
determined}: The displacement of any particle at ${\mathsf{ x}}$ is entirely
determined by ${\mathbf{ a}}\raisebox{-0.3em}{$\left|_{{\mathsf{ x}}+{{\mathbb{B}}^{(2)}}}\right.$}$.
(Heuristically
speaking, when ``deciding'' its trajectory, a particle must look not
only in a ${{\mathbb{B}}}$-neighbourhood around its current location, but
also in a ${\mathbb{B}}$-neighbourhood around each of its possible
destinations.)
\end{description}

    Several particles may be present at a given site; the rule must
assign a displacement to each of them, yielding a multiset of
displacements, which can be represented as an element of
${\mathbb{N}}^{\mathbb{B}}$.  Formally, a {\bf particle displacement
rule}\footnote{This generalizes the {\em particle automata} introduced
in \cite{MoreiraBoccaraGoles}.} (PDR) is a function
${\mathfrak{ d}}:{\mathcal{ A}}^{{\mathbb{B}}^{(2)}}{{\longrightarrow}}{\mathbb{N}}^{\mathbb{B}}$.  If ${\mathbf{ a}}\in{\mathcal{ A}}^{\mathbb{X}}$, then
for all ${\mathsf{ x}}\in{\mathbb{X}}$, we write ${\mathfrak{ d}}_{\mathsf{ x}}({\mathbf{ a}}) =
{\mathfrak{ d}}\left({\mathbf{ a}}\raisebox{-0.3em}{$\left|_{{\mathsf{ x}}+{{\mathbb{B}}^{(2)}}}\right.$}\right)$ and write the components of this object
as ${\mathfrak{ d}}_{{\mathsf{ x}}{\rightarrow}{\mathsf{ y}}}({\mathbf{ a}})$ for all ${\mathsf{ y}}\in({\mathbb{B}}+{\mathsf{ x}})$.

If $\phi\in{\mathcal{ C}}({\mathfrak{ F}};{\mathbb{N}})$, then we say that ${\mathfrak{ d}}$ is {\bf compatible} with 
$\phi$ and ${\mathfrak{ F}}$ if, for any ${\mathbf{ a}}\in{\mathcal{ A}}^{\mathbb{X}}$, with ${\mathbf{ a}}'={\mathfrak{ F}}({\mathbf{ a}})$
the following two conditions hold:
\[
\mbox{\bf(C1)} 
\quad \phi(a_{\mathsf{ o}}) \quad=\quad
\sum_{{\mathsf{ y}}\in{\mathbb{B}}} {\mathfrak{ d}}_{{\mathsf{ o}}{\rightarrow}{\mathsf{ y}}}({\mathbf{ a}})
\hspace{3em}\mbox{and}\hspace{3em}
\mbox{\bf(C2)} 
\quad
\phi(a'_{\mathsf{ o}}) \quad=\quad \sum_{{\mathsf{ x}}\in{\mathbb{B}}} {\mathfrak{ d}}_{{\mathsf{ x}}{\rightarrow}{\mathsf{ o}}}({\mathbf{ a}}).
\]
  Given $\phi$ and ${\mathfrak{ d}}$, we can
construct all $\phi$-conserving CA compatible with ${\mathfrak{ d}}$
as follows:
\begin{enumerate}
  \item Let $\displaystyle M=\max_{a\in{\mathcal{ A}}} \phi(a)$.  For every $m\in{\left[ 0..M \right]}$,
let ${\mathcal{ A}}_m=\phi^{-1}\{m\}\subset{\mathcal{ A}}$.

  \item For any ${\mathbf{ a}}\in{\mathcal{ A}}^{{\mathbb{B}}^{(2)}}$, if  $\displaystyle m=
\sum_{{\mathsf{ x}}\in{\mathbb{B}}} {\mathfrak{ d}}_{{\mathsf{ x}}{\rightarrow}{\mathsf{ o}}}({\mathbf{ a}})$ then define ${\mathfrak{ f}}({\mathbf{ a}})$ to be
any element of ${\mathcal{ A}}_m$.
\end{enumerate}
Thus, if ${{\sf card}\left[{\mathcal{ A}}_m\right]}\leq 1$ for all $m$, then the PDR uniquely
determines ${\mathfrak{ f}}$, but if ${{\sf card}\left[{\mathcal{ A}}_m\right]}> 1$ for some $m$, then there
are several CA compatible with ${\mathfrak{ d}}$.  Heuristically speaking, if
${{\sf card}\left[{\mathcal{ A}}_m\right]}\leq 1$ for all $m$, then particles have no `internal
state', and are totally interchangeable.  If ${{\sf card}\left[{\mathcal{ A}}_1\right]}>1$,
however, then a solitary particle has more than one `internal state';
\ if ${{\sf card}\left[{\mathcal{ A}}_m\right]}>1$, then a pile of $m$ particles at a site can be
arranged in more than one `configuration'.  The PDR \ ${\mathfrak{ d}}$ assigns
locations to the particles, but the local rule ${\mathfrak{ f}}$ is still
responsible for deciding what `configurations' and `internal states'
they assume.

  Note that, in general, this method yields a local rule
${\mathfrak{ f}}:{\mathcal{ A}}^{{\mathbb{B}}^{(2)}}{{\longrightarrow}}{\mathcal{ A}}$.  However, if ${\mathfrak{ d}}$ is compatible with a
cellular automaton ${\widetilde{\mathfrak{ F}}}$ having local rule
${\widetilde{\mathfrak{ f}}}:{\mathcal{ A}}^{\mathbb{B}}{{\longrightarrow}}{\mathcal{ A}}$, then condition {\bf (C2)} ensures that $\displaystyle
\sum_{{\mathsf{ x}}\in{\mathbb{B}}} {\mathfrak{ d}}_{{\mathsf{ x}}{\rightarrow}{\mathsf{ o}}}({\mathbf{ a}})$ will also depend only
on ${\mathbf{ a}}\raisebox{-0.3em}{$\left|_{{\mathbb{B}}}\right.$}$; \ thus, ${\mathfrak{ f}}$ will be a function
${\mathfrak{ f}}:{\mathcal{ A}}^{\mathbb{B}}{{\longrightarrow}}{\mathcal{ A}}$.

  The problem of constructing a PPCA is thus reduced to the problem of
constructing a PDR.  To show that {\em every} PPCA arises in this
manner, it suffices to show that every PPCA has a compatible PDR.
For arbitrary ${\mathbb{X}}$, this is surprisingly difficult; \ there are
many directions a particle can go in, and potentially several
particles vying for each destination.  When ${\mathbb{X}}={\mathbb{Z}}$, 
the one-dimensional topology obviates these complications.

\begin{prop}
\label{Prop:Z.PDR}  
If ${\mathfrak{ F}}:{\mathcal{ A}}^{\mathbb{Z}}\!\longrightarrow{\mathcal{ A}}^{\mathbb{Z}}\!$ is a cellular automaton, and
 $\phi\in{\mathcal{ C}}({\mathfrak{ F}};{\mathbb{N}})$, then there is a PDR compatible with
${\mathfrak{ F}}$ and $\phi$\hrulefill\ensuremath{\Box}
 \end{prop}

\medskip

 We will construct the PDR via a naturally defined
``flux'' function, which describes the flow rate of particles
past each point in ${\mathbb{Z}}$.
Assume ${\mathbb{B}}={\left[ -B...B \right]}$.  If ${\mathsf{ z}}\in{\mathbb{Z}}$,
then we define the {\bf flux} from ${\mathsf{ z}}$
to ${\mathsf{ z}}+1$ as follows.  Let ${\mathbf{ a}}\in{\mathcal{ A}}^{<{\mathbb{X}}}$ and let
${\mathbf{ a}}'={\mathfrak{ F}}({\mathbf{ a}})$.  Since $\phi$ is conserved, we know that
\[
  \sum_{{\mathsf{ y}}=-{\infty}}^{{\mathsf{ z}}} \phi(a_{\mathsf{ y}}) +
\sum_{{\mathsf{ y}}={\mathsf{ z}}+1}^{{\infty}} \phi(a_{\mathsf{ y}})
\ = \ \sum_{{\mathsf{ y}}=-{\infty}}^{{\infty} } \phi(a_{\mathsf{ y}})
\ = \ \sum_{{\mathsf{ y}}=-{\infty}}^{{\infty} } \phi(a'_{\mathsf{ y}})
\ = \ 
\sum_{{\mathsf{ y}}=-{\infty}}^{{\mathsf{ z}} } \phi(a'_{\mathsf{ y}}) + 
\sum_{{\mathsf{ y}}={\mathsf{ z}}+1}^{\infty} \phi(a'_{\mathsf{ y}})
\] 
  Thus,
\begin{equation}
\label{free.flux.defn}
 \sum_{{\mathsf{ y}}=-{\infty}}^{{\mathsf{ z}} } \phi(a_{\mathsf{ y}})
- 
\sum_{{\mathsf{ y}}=-{\infty}}^{{\mathsf{ z}} } \phi(a'_{\mathsf{ y}})
\ = \ 
\sum_{{\mathsf{ y}}={\mathsf{ z}}+1}^{{\infty}  } \phi(a'_{\mathsf{ y}})
- 
\sum_{{\mathsf{ y}}={\mathsf{ z}}+1}^{{\infty}  } \phi(a_{\mathsf{ y}})
\end{equation}
Let ${\vec{{\mathbf{ I}}}_{{\mathsf{ z}}}\left({\mathbf{ a}}\right)}$ be the quantity on either side of
(\ref{free.flux.defn}); this is the flow from ${\mathsf{ z}}$ to ${\mathsf{ z}}+1$.
Thus, ${\stackrel{\leftharpoonup}{{\mathbf{ I}}}_{{\mathsf{ z}}}\left({\mathbf{ a}}\right)} = - {\vec{{\mathbf{ I}}}_{{\mathsf{ z}}-1}\left({\mathbf{ a}}\right)}$ is
the flow from ${\mathsf{ z}}$ to ${\mathsf{ z}}-1$, and
$\displaystyle{\stackrel{\leftrightarrow}{{\mathbf{ I}}}_{{\mathsf{ z}}}\left({\mathbf{ a}}\right)} =  {\stackrel{\leftharpoonup}{{\mathbf{ I}}}_{{\mathsf{ z}}}\left({\mathbf{ a}}\right)} + {\vec{{\mathbf{ I}}}_{{\mathsf{ z}}}\left({\mathbf{ a}}\right)}$ 
is the total flux  out of site ${\mathsf{ z}}$.

\begin{prop}
\label{flux.prop}   Let ${\mathbf{ a}}\in{\mathcal{ A}}^{\mathbb{M}}$ and ${\mathbf{ a}}'={\mathfrak{ F}}({\mathbf{ a}})$.  If $\phi\in{\mathcal{ C}}({\mathfrak{ F}};{\mathbb{N}})$, then:
\begin{enumerate}
\item 
${\stackrel{\leftrightarrow}{{\mathbf{ I}}}_{{\mathsf{ z}}}\left({\mathbf{ a}}\right)} \ = \ -{\partial_t\, {\phi}_{{\mathsf{ z}}}\left({\mathbf{ a}}\right)}$.

\item For any  ${\mathsf{ z}}$, the value of ${\vec{{\mathbf{ I}}}_{{\mathsf{ z}}}\left({\mathbf{ a}}\right)}$ is a function only of
 ${\mathbf{ a}}\raisebox{-0.3em}{$\left|_{{\mathbb{B}}+{\mathsf{ z}}}\right.$}$.

\item  (i) \ \ $\displaystyle {\vec{{\mathbf{ I}}}_{{\mathsf{ z}}}\left({\mathbf{ a}}\right)}\  \leq \
\sum_{{\mathsf{ y}}={\mathsf{ z}}-B}^{{\mathsf{ z}}} \phi(a_{\mathsf{ y}})$; \hspace{3em} (ii) \ 
$\displaystyle {\vec{{\mathbf{ I}}}_{{\mathsf{ z}}}\left({\mathbf{ a}}\right)}\  \leq \
\sum_{{\mathsf{ y}}={\mathsf{ z}}}^{{\mathsf{ z}}+B} \phi(a'_{\mathsf{ y}})$;

 (iii) \ $\displaystyle {\stackrel{\leftharpoonup}{{\mathbf{ I}}}_{{\mathsf{ z}}}\left({\mathbf{ a}}\right)} \ \leq \
 \sum_{{\mathsf{ y}}={\mathsf{ z}}}^{{\mathsf{ z}}+B} \phi(a_{\mathsf{ y}})$; \hspace{3em} (iv) \
 $\displaystyle {\stackrel{\leftharpoonup}{{\mathbf{ I}}}_{{\mathsf{ z}}}\left({\mathbf{ a}}\right)} \ \leq \
 \sum_{{\mathsf{ y}}={\mathsf{ z}}-B}^{{\mathsf{ z}}} \phi(a'_{\mathsf{ y}})$.

\end{enumerate}
 \end{prop}
\bprf Without loss of generality, assume ${\mathsf{ z}}=0$.

{\bf Part 1} follows from equation (\ref{free.flux.defn}) by straightforward
algebra.

{\bf Part 2}: 
Suppose ${\mathbf{ b}}\in{\mathcal{ A}}^{<{\mathbb{X}}}$ with
$\displaystyle{\mathbf{ a}}\raisebox{-0.3em}{$\left|_{{\left[ -B...B \right]}}\right.$} = {\mathbf{ b}}\raisebox{-0.3em}{$\left|_{{\left[ -B...B \right]}}\right.$}$.
Define ${\mathbf{ c}}\in{\mathcal{ A}}^{<{\mathbb{X}}}$
by $\displaystyle{\mathbf{ c}}\raisebox{-0.3em}{$\left|_{{{\left( -{\infty}...B \right]}}}\right.$}={\mathbf{ a}}\raisebox{-0.3em}{$\left|_{{{\left( -{\infty}...B \right]}}}\right.$}$
while $\displaystyle{\mathbf{ c}}\raisebox{-0.3em}{$\left|_{{{\left[ -B...{\infty} \right)}}}\right.$}={\mathbf{ b}}\raisebox{-0.3em}{$\left|_{{{\left[ -B...{\infty} \right)}}}\right.$}$.
Thus, if ${\mathbf{ b}}'={\mathfrak{ F}}({\mathbf{ b}})$ and  ${\mathbf{ c}}'={\mathfrak{ F}}({\mathbf{ c}})$, then
$\displaystyle{\mathbf{ c}}'\raisebox{-0.3em}{$\left|_{{{\left( -{\infty}...0 \right]}}}\right.$}={\mathbf{ a}}'\raisebox{-0.3em}{$\left|_{{{\left( -{\infty}...0 \right]}}}\right.$}$
while $\displaystyle{\mathbf{ c}}'\raisebox{-0.3em}{$\left|_{{{\left[ 0...{\infty} \right)}}}\right.$}={\mathbf{ b}}'\raisebox{-0.3em}{$\left|_{{{\left[ 0...{\infty} \right)}}}\right.$}$.
Thus,
\begin{eqnarray*}
{\vec{{\mathbf{ I}}}_{0}\left({\mathbf{ a}}\right)}
&=&
\sum_{{\mathsf{ y}}\leq 0  } \phi(a_{\mathsf{ y}}) - \sum_{{\mathsf{ y}}\leq 0  } \phi(a'_{\mathsf{ y}})
\quad=\quad
\sum_{{\mathsf{ y}}\leq 0  } \phi(c_{\mathsf{ y}}) - \sum_{{\mathsf{ y}}\leq 0  } \phi(c'_{\mathsf{ y}})
\quad=\quad
{\vec{{\mathbf{ I}}}_{0}\left({\mathbf{ c}}\right)}
\\&=&
\sum_{0<{\mathsf{ y}}  } \phi(c'_{\mathsf{ y}}) - \sum_{0<{\mathsf{ y}}} \phi(c_{\mathsf{ y}})
\quad=\quad
\sum_{0<{\mathsf{ y}}  } \phi(b'_{\mathsf{ y}}) - \sum_{0<{\mathsf{ y}}} \phi(b_{\mathsf{ y}})
\quad=\quad
{\vec{{\mathbf{ I}}}_{0}\left({\mathbf{ b}}\right)}.
\end{eqnarray*}

{\bf Part 3}: To prove (i)
let ${\mathbf{ b}}={\left\langle {\mathbf{ a}}\raisebox{-0.3em}{$\left|_{{\mathbb{B}}}\right.$} \right\rangle }$,
Thus, ${\vec{{\mathbf{ I}}}_{0}\left({\mathbf{ b}}\right)}={\vec{{\mathbf{ I}}}_{0}\left({\mathbf{ a}}\right)}$ by {\bf Part 2}, but
$b_{\mathsf{ y}}=\mbox{\cursive O}\in\mathbf{0}$ for all ${\mathsf{ y}}\not\in{\mathbb{B}}$, so that
$\displaystyle
{\vec{{\mathbf{ I}}}_{0}\left({\mathbf{ b}}\right)} \quad \leq \quad
 \sum_{{\mathsf{ y}}\leq 0} \phi(b_{\mathsf{ y}})
\quad =\quad \sum_{{\mathsf{ y}}=-B}^{0} \phi(b_{\mathsf{ y}})
 \quad = \quad
\sum_{{\mathsf{ y}}=-B}^{0} \phi(a_{\mathsf{ y}})$.

  Inequality (iii) is proved similarly.
For inequality (ii), suppose that ${\sf supp}\left[{\mathbf{ a}}\right]\cup{\sf supp}\left[{\mathbf{ a}}'\right]
\subset{\left[ -N...{\infty} \right]}$.  Then
\begin{eqnarray*}
 {\vec{{\mathbf{ I}}}_{0}\left({\mathbf{ a}}\right)}
&= &
\sum_{{\mathsf{ y}}\leq 0} \phi(a_{\mathsf{ y}}) - \sum_{{\mathsf{ y}}\leq 0} \phi(a'_{\mathsf{ y}})
\quad = \quad 
 \sum_{{\mathsf{ y}}=-N}^0 \phi(a_{\mathsf{ y}}) - \sum_{{\mathsf{ y}}=-N-B}^0 \phi(a'_{\mathsf{ y}})\\
&=& \sum_{{\mathsf{ y}}=-N}^0 \phi(a_{\mathsf{ y}}) - \sum_{{\mathsf{ y}}=-N-B}^{B} \phi(a'_{\mathsf{ y}})
\ \  + \ \ \sum_{{\mathsf{ y}}=1}^{B} \phi(a'_{\mathsf{ y}})\\
&\leq_{[1]}&
 \sum_{{\mathsf{ y}}=-N}^0 \phi(a_{\mathsf{ y}}) -  \sum_{{\mathsf{ y}}=-N}^0 \phi(a_{\mathsf{ y}})
 \ + \ \  \sum_{{\mathsf{ y}}=1}^{B} \phi(a'_{\mathsf{ y}})
\quad\quad = \quad\quad \sum_{{\mathsf{ y}}=1}^{B} \phi(a'_{\mathsf{ y}}),
\end{eqnarray*}
 where $[1]$ follows from Theorem \ref{inf.cond}, with
${\mathbb{W}} = {\left[ -N...0 \right]}$.  Inequality (iv) is similar.
 {\tt \hrulefill $\Box$ } \end{list}  \medskip  
 
   {\bf Part 2} of Proposition \ref{flux.prop} implies that we can well-define 
${\vec{{\mathbf{ I}}}_{{\mathsf{ z}}}\left({\mathbf{ a}}\right)}$ even when ${\mathbf{ a}}$ has infinite support.
Note that {\bf Part 2} does {\em not} imply that
\[
  {\vec{{\mathbf{ I}}}_{0}\left({\mathbf{ a}}\right)} =  \sum_{{\mathsf{ y}}=-B}^{0} \phi(a_{\mathsf{ y}})
- 
\sum_{{\mathsf{ y}}=-B}^{0} \phi(a'_{\mathsf{ y}})
\ = \ 
\sum_{{\mathsf{ y}}=1}^{B} \phi(a'_{\mathsf{ y}})
- 
\sum_{{\mathsf{ y}}=1}^{B} \phi(a_{\mathsf{ y}}). 
\]

  \medskip         \refstepcounter{thm} {\bf Example \thethm:}  \setcounter{enumi}{\thethm} \begin{list}{(\alph{enumii})}{\usecounter{enumii}} 			{\setlength{\leftmargin}{0em} 			\setlength{\rightmargin}{0em}}   
\item  Suppose ${\mathfrak{ F}}= {{{\boldsymbol{\sigma}}}^{5}} $ is the five-fold shift.  
If ${\mathbf{ a}} = \left[\ldots1011110\hat{0}110100101\ldots\right]$, then
${\mathbf{ a}} = \left[\ldots1011110011010\hat{0}101\ldots\right]$, and
${\stackrel{\leftharpoonup}{{\mathbf{ I}}}_{0}\left({\mathbf{ a}}\right)}=3$ (here, the hat indicates the
position of $a_0$).  If ${\mathbf{ a}} = \left[\ldots11111\hat{1}11111\ldots\right]$,
then ${\stackrel{\leftharpoonup}{{\mathbf{ I}}}_{0}\left({\mathbf{ a}}\right)}=5$.

\item Suppose ${\mathfrak{ F}}$ is CA \#187 from Example \ref{EXPPCA}. If ${\mathbf{ a}}
= \left[\ldots\hat{1}0\ldots\right]$, then ${\vec{{\mathbf{ I}}}_{0}\left({\mathbf{ a}}\right)}=1$, while
if ${\mathbf{ a}} = \left[\ldots\hat{1}1\ldots\right]$, then ${\vec{{\mathbf{ I}}}_{0}\left({\mathbf{ a}}\right)}=0$.
 	\hrulefill\end{list}   			

\begin{list}{} 			{\setlength{\leftmargin}{1em} 			\setlength{\rightmargin}{0em}}                         \item {\bf \hspace{-1em}  Proof of Proposition  \ref{Prop:Z.PDR}: \ \ }
  Again, assume without loss of generality that ${\mathsf{ z}}=0$.
We construct the PDR from the flux as follows:

{\bf Case 0:} ({\em ${\stackrel{\leftharpoonup}{{\mathbf{ I}}}_{0}\left({\mathbf{ a}}\right)} \leq0$ and 
 ${\vec{{\mathbf{ I}}}_{0}\left({\mathbf{ a}}\right)} \leq0$})
\quad
No particles leave site $0$ in either direction, so
${\mathfrak{ d}}_0({\mathbf{ a}})\equiv 0$.

{\bf Case 1:} ({\em $0\leq{\stackrel{\leftharpoonup}{{\mathbf{ I}}}_{0}\left({\mathbf{ a}}\right)}$ and $0\leq{\vec{{\mathbf{ I}}}_{0}\left({\mathbf{ a}}\right)}$})
\quad
Now, particles are leaving site $0$ in both directions.
For ${\mathsf{ z}}>0$, if $\displaystyle\sum_{0<{\mathsf{ y}}<{\mathsf{ z}}} \phi(a'_{\mathsf{ y}})< {\vec{{\mathbf{ I}}}_{0}\left({\mathbf{ a}}\right)}$
then define
\[ 
{\mathfrak{ d}}_{0{\rightarrow}{\mathsf{ z}}}({\mathbf{ a}}) = 
 \min\left\{ \phi(a'_{\mathsf{ z}}), \ \ {\vec{{\mathbf{ I}}}_{0}\left({\mathbf{ a}}\right)} - \sum_{0<{\mathsf{ y}}<{\mathsf{ z}}} \phi(a'_{\mathsf{ y}})\right\}
 \]
and set ${\mathfrak{ d}}_{0{\rightarrow}{\mathsf{ z}}}({\mathbf{ a}})=0$ if $\displaystyle {\vec{{\mathbf{ I}}}_{0}\left({\mathbf{ a}}\right)}\leq
\sum_{0<{\mathsf{ y}}<{\mathsf{ z}}} \phi(a'_{\mathsf{ z}})$.
Likewise, for ${\mathsf{ z}}<0$, if $\displaystyle\sum_{{\mathsf{ z}}<{\mathsf{ y}}<0} \phi(a'_{\mathsf{ z}}) <
{\stackrel{\leftharpoonup}{{\mathbf{ I}}}_{0}\left({\mathbf{ a}}\right)}$, then define
\[ 
{\mathfrak{ d}}_{0{\rightarrow}{\mathsf{ z}}}({\mathbf{ a}}) = 
 \min\left\{ \phi(a'_{\mathsf{ z}}), \ \ {\stackrel{\leftharpoonup}{{\mathbf{ I}}}_{0}\left({\mathbf{ a}}\right)}  - \sum_{{\mathsf{ z}}<{\mathsf{ y}}<0} \phi(a'_{\mathsf{ y}})\right\}
\]
 and set ${\mathfrak{ d}}_{0{\rightarrow}{\mathsf{ z}}}({\mathbf{ a}})=0$ if 
$\displaystyle{\stackrel{\leftharpoonup}{{\mathbf{ I}}}_{0}\left({\mathbf{ a}}\right)} \leq \sum_{{\mathsf{ z}}<{\mathsf{ y}}<0} \phi(a'_{\mathsf{ z}})$.

{\bf Case 2.1:} (${\stackrel{\leftharpoonup}{{\mathbf{ I}}}_{0}\left({\mathbf{ a}}\right)}\leq 0<{\vec{{\mathbf{ I}}}_{0}\left({\mathbf{ a}}\right)}$)
\quad
Now, particles enter $0$ from the left, and leave to the right.
Let $J_0=\phi(a'_0)$; the first $J_0$ particles entering $0$ from
the left will fill the $J_0$ available destinations
at $0$.  Let $J_1=\max\{0, \ -{\stackrel{\leftharpoonup}{{\mathbf{ I}}}_{0}\left({\mathbf{ a}}\right)}-J_0\}$; \ the
next $J_1$ particles entering from the left must pass through $0$,
and will fill up the next $J_1$ destinations available to the right of
$0$.

${\vec{{\mathbf{ I}}}_{0}\left({\mathbf{ a}}\right)}$ particles leave $0$ to the right.
If $J_1>0$, then the first $J_1$ of these
particles are from ${\left( -{\infty}...0 \right)}$, 
while the last ${\vec{{\mathbf{ I}}}_{0}\left({\mathbf{ a}}\right)}-J_1$ are the particles
originating at $0$ itself ({\bf Part 1} of Proposition
\ref{flux.prop} implies that ${\vec{{\mathbf{ I}}}_{0}\left({\mathbf{ a}}\right)}-J_1 \ = \ {\vec{{\mathbf{ I}}}_{0}\left({\mathbf{ a}}\right)}
+{\stackrel{\leftharpoonup}{{\mathbf{ I}}}_{0}\left({\mathbf{ a}}\right)} + \phi(a'_0) \ = \ -{\partial_t\, {\phi}_{0}\left({\mathbf{ a}}\right)} +
\phi(a'_0) \ =  \ \phi(a_0)$).

   Thus, for all ${\mathsf{ z}}>0$, if
$\displaystyle J_1< \sum_{0<{\mathsf{ y}}<{\mathsf{ z}}} \phi(a'_{\mathsf{ y}}) <{\vec{{\mathbf{ I}}}_{0}\left({\mathbf{ a}}\right)}$
then define
\[ 
{\mathfrak{ d}}_{0{\rightarrow}{\mathsf{ z}}}({\mathbf{ a}}) = 
 \min\left\{ \phi(a'_{\mathsf{ z}}), \ \ {\vec{{\mathbf{ I}}}_{0}\left({\mathbf{ a}}\right)} - \sum_{0<{\mathsf{ y}}<{\mathsf{ z}}} \phi(a'_{\mathsf{ y}})\right\}
 \]
and set ${\mathfrak{ d}}_{0{\rightarrow}{\mathsf{ z}}}({\mathbf{ a}})=0$ if 
$\displaystyle\sum_{0<{\mathsf{ y}}\leq {\mathsf{ z}}} \phi(a'_{\mathsf{ y}})\leq J_1$ or if
$\displaystyle {\vec{{\mathbf{ I}}}_{0}\left({\mathbf{ a}}\right)} \leq \sum_{0<{\mathsf{ y}}<{\mathsf{ z}}} \phi(a'_{\mathsf{ y}})$.  This leaves the
boundary case when 
$\displaystyle\sum_{0<{\mathsf{ y}}< {\mathsf{ z}}} \phi(a'_{\mathsf{ y}})< J_1
< \sum_{0<{\mathsf{ y}}\leq{\mathsf{ z}}} \phi(a'_{\mathsf{ y}})$.  In this case, let
$\displaystyle
{\mathfrak{ d}}_{0{\rightarrow}{\mathsf{ z}}}({\mathbf{ a}}) = \sum_{0<{\mathsf{ y}}\leq{\mathsf{ z}}} \phi(a'_{\mathsf{ y}}) - J_1$.

{\bf Case 2.2:} (${\vec{{\mathbf{ I}}}_{0}\left({\mathbf{ a}}\right)}\leq0<{\stackrel{\leftharpoonup}{{\mathbf{ I}}}_{0}\left({\mathbf{ a}}\right)}$)
\quad
Now, particles enter $0$ from the right, and leave to the
left.  This case is handled analogously to Case 2.1

\medskip

  It remains to verify that ${\mathfrak{ d}}$ is a PDR, and is compatible
with ${\mathfrak{ F}}$.  
In each of Cases 1, 2.1, and 2.2, ${\mathfrak{ d}}_{0{\rightarrow}{\mathsf{ z}}}({\mathbf{ a}})=0$ for any
${\mathsf{ z}}\not\in{\mathbb{B}}$; this follows from inequalities (ii) and (iv) in
{\bf Part 3} of Proposition \ref{flux.prop}.  Thus, condition {\bf
(D2)} is satisfied.  Also, the value of ${\mathfrak{ d}}_{0}({\mathbf{ a}})$ is
determined by ${\mathbf{ a}}'\raisebox{-0.3em}{$\left|_{{\mathbb{B}}}\right.$}$, ${\vec{{\mathbf{ I}}}_{0}\left({\mathbf{ a}}\right)}$ and
${\stackrel{\leftharpoonup}{{\mathbf{ I}}}_{0}\left({\mathbf{ a}}\right)}$.  Clearly, ${\mathbf{ a}}'\raisebox{-0.3em}{$\left|_{{\mathbb{B}}}\right.$}$ is determined by
${\mathbf{ a}}\raisebox{-0.3em}{$\left|_{{\mathbb{B}}^{(2)}}\right.$}$, while
{\bf Part 2} of Proposition \ref{flux.prop} says that
${\vec{{\mathbf{ I}}}_{0}\left({\mathbf{ a}}\right)}$ and ${\stackrel{\leftharpoonup}{{\mathbf{ I}}}_{0}\left({\mathbf{ a}}\right)}$ are determined by
${\mathbf{ a}}\raisebox{-0.3em}{$\left|_{{\mathbb{B}}}\right.$}$.  Thus, ${\mathfrak{ d}}_{0}({\mathbf{ a}})$
is a function only of ${\mathbf{ a}}\raisebox{-0.3em}{$\left|_{{\mathbb{B}}^{(2)}}\right.$}$, so  
condition {\bf (D3)} holds.
We have described the algorithm at ${\mathsf{ z}}=0$, but
 we apply the same algorithm at all points, so condition
{\bf (D1)} holds automatically.

 Compatibility condition {\bf (C1)} holds by construction in each Case. 
To check {\bf (C2)} suppose ${\mathbf{ a}}\in{\mathcal{ A}}^{<{\mathbb{Z}}}$, with
${\sf supp}\left[{\mathbf{ a}}\right]={\left[ -N...N \right]}$.  Beginning at $-N$ and proceeding to the
right, we can inductively verify {\bf (C2)}
for each ${\mathsf{ z}}\in{\mathbb{Z}}$ by applying {\bf Part 1} of Proposition
\ref{flux.prop}.
 {\tt \hrulefill $\Box$ } \end{list}  \medskip  

\section{Conclusion}

  In \S\ref{S:finitary} and \S\ref{S:measure},
we provide practical methods for detecting the existence of
conservation laws, while in \S\ref{S:infinitary} and \S\ref{S:Cesaro},
we provide abstract characterizations of such laws. 
Theorem \ref{inf.cond} from \S\ref{S:infinitary} is used in 
\S\ref{S:construct} to develop a notion of `flux', which
yields a method for constructing cellular automata having a particular
conservation law. 

  However, the `displacement representation' constructed in
\S\ref{S:construct} is inapplicable to the case ${\mathbb{X}}={\mathbb{Z}}^D$,
$D\geq 2$.  Do displacement representations exist for PPCA on higher
dimensional lattices?  If we interpret a conserved quantity as the
density of some material, many questions remain about the
`hydrodynamics' of this material: \ its patterns of flow,
concentration, and diffusion.  We also expect that higher-dimensional
PPCA may exhibit complex particle dynamics, including the formation of
complex, large-scale, stable clusters analogous to molecules.  What is
a good framework for studying these quasichemical dynamics?

\nocite{Moreira}
{\small
\bibliographystyle{plain}
\bibliography{bibliography}

\begin{thebibliography}{10}

\bibitem{MoreiraBoccaraGoles}
{Nino Boccara} {Andr\'es Moreira} and {Eric Goles}.
\newblock Number-conserving one-dimensional cellular automata and particle
  representation.
\newblock {\em {\rm (preprint)}}, January 2002.

\bibitem{DurandFormentiRoka}
{Enrico Formenti} {Bruno Durand} and {Zsuzsanna R\'oka}.
\newblock Number conserving cellular automata: from decidability to dynamics.
\newblock {\em {\rm preprint available at} {\tt
  http://arXiv.org/ps/nlin.CG/0102035}}, 2001.

\bibitem{HedlundCA}
G.~Hedlund.
\newblock Endomorphisms and automorphisms of the shift dynamical systems.
\newblock {\em Mathematical System Theory}, 3:320--375, 1969.

\bibitem{EsserSchreckenberg}
{J. Esser} and {M. Schreckenberg}.
\newblock Microscopic simulations of urban traffic based on cellular automata.
\newblock {\em International Journal of Modern Physics C}, 8:1025--1036, 1997.

\bibitem{NagelSchrekenberg}
{K. Nagel} and {M. Schreckenberg}.
\newblock A cellular automaton model for freeway traffic.
\newblock {\em J. Physique I}, 2:2221, 1992.

\bibitem{Kohyama1}
Tamotsu Kohyama.
\newblock Cellular automata with particle conservation.
\newblock {\em Progress of Theoretical Physics}, 81(1):47--59, January 1989.

\bibitem{Kohyama2}
Tamotsu Kohyama.
\newblock Cluster growth in particle-conserving cellular automata.
\newblock {\em Journal of Statistical Physics}, 63(3/4):637--651, 1991.

\bibitem{KotzeSteeb}
{L. Kotze} and {W. H. Steeb}.
\newblock Conservation laws in cellular automata.
\newblock In {P. G. L. Leach} and {W.H. Steeb}, editors, {\em Finite
  Dimensional Integrable Nonlinear Dynamical Systems}, pages 333--346, New
  Jersey, 1988. World Scientific.

\bibitem{FukuiIshibashi}
{M. Fukui} and {Y. Ishibashi}.
\newblock Traffic flow in {1D} cellular automaton model including cars moving
  with high speed.
\newblock {\em J. Phys. Soc. Japan}, 65:1868--1870, 1996.

\bibitem{Moreira}
Andr{\'e}s Moreira.
\newblock Universality and decidability of number-conserving cellular automata.
\newblock {\em {\rm submitted to} Theoretical Computer Science}, August 2001.

\bibitem{BoccaraFuks1}
{Nino Boccara} and {Henryk Fuk{\'s}}.
\newblock Cellular automaton rules conserving the number of active sites.
\newblock {\em Journal of Physics A: Math. Gen.}, 31:6007--6018, 1998.

\bibitem{BoccaraFuks2}
{Nino Boccara} and {Henryk Fuk{\'s}}.
\newblock Number-conserving cellular automaton rules.
\newblock {\em Fundamenta Informaticae}, pages 1--14, 2000.

\bibitem{SimonNagel}
{P. M. Simon} and {K. Nagel}.
\newblock Simplified cellular automaton model for city traffic.
\newblock {\em Physical Review E}, 58:1286--1295, 1998.

\bibitem{PivatoCJM}
Marcus Pivato.
\newblock Building a stationary stochastic process from a finite-dimensional
  marginal.
\newblock {\em Canadian Journal of Mathematics}, 53(2):382--413, 2001.

\bibitem{TakesueERCA}
Shinji Takesue.
\newblock Ergodic properties and thermodynamic behaviour of elementary
  reversible cellular automata.
\newblock {\em Journal of Statistical Physics}, 56:371, 1989.

\bibitem{Tempelman}
A.~A. Tempel'man.
\newblock Ergodic theorems for general dynamical systems.
\newblock {\em Soviet Math. Doklady}, 8(5):1213--1216, 1967.
\newblock (English Translation).

\bibitem{HattoriTakesue}
{Tetsuya Hattori} and {Shinji Takesue}.
\newblock Additive conserved quantities in discrete-time lattice dynamical
  systems.
\newblock {\em Physica D}, 49:295--322, 1991.

\bibitem{WolframBook}
Stephen Wolfram.
\newblock {\em Cellular Automata and Complexity}.
\newblock Addison-Wesley, Reading, Massachusetts, 1994.

\end{thebibliography}
}
\end{document}